\documentclass[12pt]{article}

\usepackage{amsmath}
\usepackage{amsfonts}
\usepackage{amssymb}
\usepackage{graphicx}

\makeatletter
\@addtoreset{equation}{section}

\makeatother
\evensidemargin\oddsidemargin
\newtheorem{prop}{Proposition}

\newtheorem{theo}{Theorem}
\newtheorem{theor}{Theorem}[section]
\newtheorem{fact}{Fact}

\begin{document}

\title{\textbf{Uniform in bandwidth consistency of conditional $U$-statistics}}
\author{J. Dony* and D.M. Mason${}^1$}
\date{\vspace{-2cm}  }
\maketitle
\begin{align*}
\hspace{0.5cm}&\textit{{\small ${}^*$Mathematics Department, Free University of Brussels (VUB). Pleinlaan 2,}}\\
 &\textit{{B-1050 Brussels, Belgium.  e-mail : jdony@vub.ac.be}}\\
& \textit{{\small  ${}^1$Food and Resource Economics, University of Delaware. 206 Townsend Hall,}}\\
&\textit{{ Newark, DE 19717.  e-mail : davidm@udel.edu}} 
\end{align*}

\begin{abstract}
{\footnotesize In 1991 Stute introduced a class of estimators called
conditional $U$--statistics. They can be seen as a generalization of the Nadaraya-Watson estimator for the regression function, and he proved their strong pointwise consistency
to
\begin{equation*}
m(\mathbf{t}):=\mathbb{E}[ g(Y_{1},\ldots,Y_{m})|(X_{1},\ldots, X_{m})=%
\mathbf{t}],\quad\mathbf{t}\in\mathbb{R}^{m}.
\end{equation*}
Very recently, Gin\'{e} and Mason introduced the notion of a
local $U $--process, which generalizes that of a local empirical process,
and obtained central limit theorems and laws of the iterated logarithm for
this class. We apply the methods developed in Einmahl and Mason (2005)
and Gin\'{e} and Mason (2007a,b) to establish uniform in bandwidth
consistency to $m(\mathbf{t})$ of the estimator proposed by Stute. }
\end{abstract}


\section{Introduction and statement of main results}

Let $(X,Y)$, $(X_{1},Y_{1}),\ldots,(X_{n},Y_{n})$ be independent random
vectors with common joint density function $f:\mathbb{R}\times\mathbb{R}%
\rightarrow[ 0,\infty[ $, and for a measurable function $\varphi:%
\mathbb{R}^{m}\rightarrow\mathbb{R}$, consider the regression function 
\begin{equation*}
m_{\varphi}(\mathbf{t})=\mathbb{E}\left[ \varphi(Y_{1},\dots,Y_{m})|(X_{1},%
\dots,X_{m})=\mathbf{t}\right] ,\quad\mathbf{t}\in\mathbb{R}^{m}.
\end{equation*}
Stute \cite{STUT91} introduced a class of estimators for $m_{\varphi }(%
\mathbf{t})$, called conditional $U$--statistics and defined pointwise for $%
\mathbf{t}\in\mathbb{R}^{m}$ as 
\begin{equation}
\widehat{m}_{n}(\mathbf{t};h_{n})=\frac{\sum_{(i_{1},\dots,i_{m})\in
I_{n}^{m}}\varphi(Y_{i_{1}},\dots,Y_{i_{m}})K\left( \frac{t_{1}-X_{i_{1}}}{%
h_{n}}\right) \cdots K\left( \frac{t_{m}-X_{i_{m}}}{h_{n}}\right) }{\sum
_{(i_{1},\dots,i_{m})\in I_{n}^{m}}K\left( \frac{t_{1}-X_{i_{1}}}{h_{n}}%
\right) \cdots K\left( \frac{t_{m}-X_{i_{m}}}{h_{n}}\right) },
\label{def:estimator}
\end{equation}
where 
\begin{equation}
I_{n}^{m}=\{(i_{1},\ldots,i_{m}):1\leq i_{j}\leq n,i_{j}\neq i_{l}\;\text{if}%
\;j\neq l\},  \label{II}
\end{equation}
and $0<h_{n}<1$ goes to zero at a certain rate. Soon afterwards, Sen \cite%
{SEN94} obtained results on uniform consistency of this estimator. We shall
adapt and extend the methods developed in Einmahl and Mason \cite{EM2005}
and Gin\'{e} and Mason \cite{GM07a,GM07b} to show that under appropriate
regularity conditions a much stronger form of consistency holds, namely
uniform in bandwidth consistency of $\widehat{m}_{n}$. This means that with
probability $1$, 
\begin{equation}
\limsup_{n\rightarrow\infty}\sup_{\widetilde a_{n}\leq h\leq b_{n}}\sup_{%
\mathbf{t}\in\lbrack c,d]^{m}}|\widehat{m}_{n}(\mathbf{t};h)-m_{\varphi }(%
\mathbf{t})|=0,  \label{M}
\end{equation}
for $-\infty<c<d<\infty$ and $\widetilde a_{n}<b_{n}$, as long as $%
\widetilde a_{n}\rightarrow0$, $b_{n}\rightarrow0$ and $b_{n}/\widetilde
a_{n}\rightarrow\infty$ at rates depending upon the moments of $%
\varphi(Y_{1},\dots,Y_{m})$. Moreover, we shall show that (\ref{M}) holds
uniformly in $\varphi\in\mathcal{F}$ as well. In fact, our results extend
those of Einmahl and Mason \cite{EM2005}, who treat the case $m=1$.\newline

We shall infer (\ref{M}) via general uniform in bandwidth results for a
specific $U-$statistic process indexed by a class of functions. We define
this process in (\ref{LU}) below. Towards this end, for $m\leq n$, consider
a class $\mathcal{F}$ of measurable functions $g:\mathbb{R}^{m}\rightarrow\mathbb{R}$ such that $\mathbb{E} g^{2}(Y_{1},\dots,Y_{m})<\infty$, which satisfies the following conditions $(F.i)-(F.iii)$. First, to avoid measurability problems, we assume that 
\begin{equation*}
\mathcal{F}\;\text{ is a pointwise measurable class,}\leqno{(F.i)}
\end{equation*}
i.e. there exists a countable subclass $\mathcal{F}_{0}$ of $\mathcal{F}$
such that we can find for any function $g\in\mathcal{F}$ a sequence of
functions $g_{m}\in\mathcal{F}_{0}$ for which $g_{m}(z)\rightarrow g(z),z\in 
\mathbb{R}^{m}$. This condition is discussed in van der Vaart and Wellner 
\cite{VVW96}. We also assume that $\mathcal{F}$ has a measurable envelope
function %\be\label{FF}%
\begin{equation*}
F(\mathbf{y})\geq\sup_{g\in\mathcal{F}}|g(\mathbf{y})|,\quad\mathbf{y}\in%
\mathbb{R}^{m}.\leqno{(F.ii)}
\end{equation*}
%\end{equation}
Finally we assume that $\mathcal{F}$ is of VC--type with characteristics $A$
and $v$ (\textquotedblleft VC" for Vapnik and \v{C}ervonenkis), meaning that
for some $A\geq3$ and $v\geq1$, %\be\label{NN}%
\begin{equation*}
\mathcal{N}(\mathcal{F},L_{2}(Q),\varepsilon)\leq\left( \frac{A\Vert
F\Vert_{L_{2}(Q)}}{\varepsilon}\right) ^{v},\quad0<\varepsilon\leq2\Vert
F\Vert_{L_{2}(Q)},\leqno{(F.iii) }
\end{equation*}
%\end{equation}
where for $\varepsilon>0$, $\mathcal{N}(\mathcal{F},L_{2}(Q),\varepsilon)$
is defined as the smallest number of $L_{2}(Q)$ open balls of radius $%
\varepsilon$ required to cover $\mathcal{F}$, and $Q$ is any probability
measure on $(\mathbb{R}^{m},\mathcal{B})$ such that $\Vert
F\Vert_{L_{2}(Q)}<\infty$. (If $(F.iii)$ holds for $\mathcal{F}$, then we
say that the VC--type class $\mathcal{F}$ admits the characteristics $A$ and 
$v$.)\newline

Let now $K:\mathbb{R}\rightarrow\mathbb{R}$ be a kernel function with
support contained in $[-1/2,1/2]$ satisfying 
\begin{equation*}
\sup_{x\in\mathbb{R}}|K(x)|=:\kappa<\infty\quad\text{and}\quad\int K(x)dx=1.
\leqno{(K.i)}
\end{equation*}
For such kernels, we consider the class of functions  $\mathcal{K}:=\{hK_{h}(t-\cdot):h>0,t\in\mathbb{R}\}$ and assume that 
\begin{equation*}
\mathcal{K} \text{ is pointwise measurable and of VC--type }, \leqno{(K.ii)}
\end{equation*}

\noindent where as usual $K_{h}(z)=h^{-1}K(z/h)$, $z\in \mathbb{R}$. Furthermore, let 
\begin{equation*}
\widetilde{K}(\mathbf{t}):=\prod_{j=1}^{m}K(t_{j})\leqno{(K.iii)}
\end{equation*}
denote the product kernel. Next, if $(S,\mathcal{S})$ is a measurable space, define the general $U$--statistic with
kernel $H:S^{k}\rightarrow \mathbb{R}$ based on $S$--valued random variables 
$Z_{1},\ldots ,Z_{n}$ as 
\begin{equation*}
U_{n}^{(k)}(H):=\frac{(n-k)!}{n!}\sum_{\mathbf{i}\in
I_{n}^{k}}H(Z_{i_{1}},\ldots ,Z_{i_{k}}),\quad 1\leq k\leq n,
\end{equation*}%
where $I_{n}^{k}$ is defined as in (\ref{II}) with $m=k$. (Note that we do
not require $H$ to be symmetric here.) For a bandwidth $0<h<1$ and $g\in 
\mathcal{F}$, consider the $U$--kernel 
\begin{equation*}
G_{g,h,\mathbf{t}}(\mathbf{x},\mathbf{y}):=g(\mathbf{y})\widetilde{K}_{h}(%
\mathbf{t}-\mathbf{x}),\quad \mathbf{x},\mathbf{y},\mathbf{t}\in \mathbb{R}%
^{m},
\end{equation*}%
and for the sample $(X_{1},Y_{1}),\ldots ,(X_{n},Y_{n})$, define 
\begin{equation*}
U_{n}(g,h,\mathbf{t}):=U_{n}^{(m)}(G_{g,h,\mathbf{t}})=\frac{(n-m)!}{n!}%
\sum_{\mathbf{i}\in I_{n}^{m}}G_{g,h,\mathbf{t}}(\mathbf{X}_{\mathbf{i}},%
\mathbf{Y}_{\mathbf{i}}),
\end{equation*}%
where throughout this paper we shall use the notation 
\begin{align*}
& \mathbf{X}=(X_{1},\ldots ,X_{m})\in\mathbb{R}^m\quad \text{ and }\quad \mathbf{X}_{%
\mathbf{i}}:=(X_{i_{1}},\ldots ,X_{i_{k}})\in\mathbb{R}^k,\quad \mathbf{i}\in I_{n}^{k}, \\
& \mathbf{Y}=(Y_{1},\ldots ,Y_{m})\in\mathbb{R}^m\quad \text{ and }\quad \mathbf{Y}_{%
\mathbf{i}}:=(Y_{i_{1}},\ldots ,Y_{i_{k}})\in\mathbb{R}^k,\quad \mathbf{i}\in I_{n}^{k}.
\end{align*}%
Now introduce the \textit{$U$--statistic process} 
\begin{equation}
u_{n}(g,h,\mathbf{t}):=\sqrt{n}\{U_{n}(g,h,\mathbf{t})-\mathbb{E}U_{n}(g,h,%
\mathbf{t})\}.  \label{LU}
\end{equation}

We shall establish a strong uniform in bandwidth consistency result for the $%
U$--statistic process in (\ref{LU}). Theorem \ref{theo:bounded} gives such a
result for bounded classes of functions $\mathcal{F}$, while Theorem \ref%
{theo:unbounded} is applicable for unbounded classes $\mathcal{F}$ which
satisfy a conditional moment condition stated in (\ref{unbounded}) below. In
the bounded case, we assume that the envelope function of $\mathcal{F}$ is
bounded by some finite constant $M$, i.e., (\ref{bounded}) holds.

\begin{theo}
\label{theo:bounded} Suppose that the marginal density $f_{X}$ of $X$ is
bounded, and  let $a_{n}=c(\log n/n)^{1/m}$ for $c>0$. If the class of functions $\mathcal{F}$ is bounded in the sense
that for some $0<M<\infty$, 
\begin{equation}
F(\mathbf{y})\leq M,\quad\mathbf{y}\in\mathbb{R}^{m},  \label{bounded}
\end{equation}
we can infer under the above mentioned assumptions on $\mathcal{F}$ and $%
\mathcal{K}$ that for all $c>0$ and $0<b_{0}<1$ there exists a constant $%
0<C<\infty$ such that 
\begin{equation*}
\limsup_{n\rightarrow\infty}\sup_{a_{n}\leq h\leq b_{0}}\sup_{g\in\mathcal{F}%
}\sup_{\mathbf{t}\in\mathbb{R}^{m}}\frac{\sqrt{nh^{m}}|U_{n}(g,h,\mathbf{t})-%
\mathbb{E}U_{n}(g,h,\mathbf{t})|}{\sqrt{|\log h|\vee\log\log n}}\leq C,\quad%
\text{a.s.}
\end{equation*}
\end{theo}
\bigskip
\begin{theo}
\label{theo:unbounded} Suppose that the marginal density $f_{X}$ of $X$ is
bounded, and for $c>0$ let $a_{n}^{\prime}=c((\log n/n)^{1-2/p})^{1/m}$. If $\mathcal{F}$ is unbounded but satisfies for some $p>2$ 
\begin{equation}
\mu_{p}:=\sup_{\mathbf{x}\in\mathbb{R}^{m}}\mathbb{E}[F^{p}(\mathbf{Y})|\:\mathbf{X}=\mathbf{x}]<\infty,
\label{unbounded}
\end{equation}
we can infer under the above mentioned assumptions on $\mathcal{F}$ and $%
\mathcal{K}$ that for all $c>0$ and $0<b_{0}<1$ there exists a constant $%
0<C^{\prime}<\infty$ such that, 
\begin{equation*}
\limsup_{n\rightarrow\infty}\sup_{a_{n}^{\prime}\leq h\leq b_{0}}\sup _{g\in%
\mathcal{F}}\sup_{\mathbf{t}\in\mathbb{R}^{m}}\frac{\sqrt{nh^{m}}|U_{n}(g,h,%
\mathbf{t})-\mathbb{E}U_{n}(g,h,\mathbf{t})|}{\sqrt{|\log h|\vee\log\log n}}%
\leq C^{\prime},\quad\text{a.s.}
\end{equation*}
\end{theo}

\bigskip

From now on, we shall write $\hat{m}_{n,\varphi}(\mathbf{t},h)$ for the
estimator of the regression function defined in (\ref{def:estimator}) to
stress the role of $\varphi(\mathbf{y})$. It is clear that $\hat{m}_{n,\varphi }(%
\mathbf{t},h)$ can be rewritten for all $\varphi\in\mathcal{F}$ as 
\begin{equation*}
\hat{m}_{n, \varphi}(\mathbf{t},h)=\frac{\sum_{\mathbf{i}\in
I_{n}^{m}}\varphi(\mathbf{Y}_{\mathbf{i}})\widetilde{K}_{h}(\mathbf{t}-%
\mathbf{X}_{\mathbf{i}})}{\sum_{\mathbf{i}\in I_{n}^{m}}\widetilde{K}_{h}(%
\mathbf{t}-\mathbf{X}_{\mathbf{i}})}=\frac{U_{n}(\varphi,h,\mathbf{t})}{%
U_{n}(1,h,\mathbf{t})},
\end{equation*}
where we denote by $U_{n}(1,h,\mathbf{t})$ the $U$--statistic $U_{n}(g,h,%
\mathbf{t})$ with $g\equiv1$. To prove the uniform consistency of $\hat{m}%
_{n, \varphi}(\mathbf{t},h)$ to $m_{\varphi}(\mathbf{t})$, we shall consider
another, but more appropriate, centering factor than the expectation $%
\mathbb{E}\hat{m}_{n, \varphi}(\mathbf{t},h)$, which may not exist or be
difficult to compute. Define the centering 
\begin{equation}
\widehat{\mathbb{E}}\hat{m}_{n, \varphi}(\mathbf{t},h):=\frac{\mathbb{E}%
U_{n}(\varphi,h,\mathbf{t})}{\mathbb{E}U_{n}(1,h,\mathbf{t})}.
\label{def:cent}
\end{equation}
This centering permits us to apply Theorems 1 and 2 (depending on whether
the class $\mathcal{F}$ is bounded in the sense of (\ref{bounded}) or
unbounded in the sense of (\ref{unbounded})) to derive results on the
convergence rates of the process $\hat{m}_{n, \varphi}(\mathbf{t},h)-%
\widehat{\mathbb{E}}\hat
{m}_{n, \varphi}(\mathbf{t},h)$ to zero and the
consistency of $\hat{m}_{n, \varphi}(\mathbf{t},h)$, uniformly in bandwidth.%
\newline

For any compact interval $I=[ c,d] $ with $-\infty<c<d<\infty$ and $\eta>0$,
define $I^{\eta}=\left[ c-\eta,d+\eta\right] $ and denote as usual the
marginal density function of $X$ by $f_{X}$. Then introduce the class of
functions defined on the compact subset $J^{m}=I^{\eta}\times\ldots\times
I^{\eta}$ of $\mathbb{R}^{m}$, 
\begin{equation}
\mathcal{M}=\{m_{\varphi}(\cdot)\widetilde{f}(\cdot):\varphi\in\mathcal{F}\},
\label{equi}
\end{equation}
where the function $\widetilde{f}$ $:\mathbb{R}^{m}\rightarrow\mathbb{R}$ is
defined as 
\begin{equation}
\widetilde{f}\left( \mathbf{t}\right) :=\int f(t_{1},y_{1})\ldots
f(t_{m},y_{m})dy_{1}\ldots dy_{m}=f_{X}(t_{1})\ldots f_{X}(t_{m}).
\label{ff}
\end{equation}

\noindent We have now introduced all the notation that we need to state our
results on the uniform consistency of the conditional $U-$statistic
estimator proposed by Stute for the general regression function, where this
consistency is uniform in $\varphi\in\mathcal{F}$ and in bandwidth as well.

\begin{theo}
\label{theo:cons} Besides being bounded, suppose that the marginal density
function $f_{X}$ of $X$ is continuous and strictly positive on the interval $%
J=I^{\eta}$, where $I=[c,d] $ is a compact interval and $\eta>0$. Assume
that the class of functions $\mathcal{M}$ is uniformly equicontinuous. Then
it follows that for all sequences $0<b_{n}<1$ with $b_{n}\rightarrow0$, 
\begin{equation*}
\sup_{0< h<b_{n}}\sup_{\varphi\in\mathcal{F}}\sup_{\mathbf{t}\in I^{m}}|%
\widehat{\mathbb{E}}\hat{m}_{n, \varphi}(\mathbf{t},h)-m_{\varphi }(\mathbf{t%
})|=o(1),
\end{equation*}
where $I^{m}=I\times\ldots\times I$.
\end{theo}

\begin{theo}
\label{theo:conv-rates} Besides being bounded, suppose that the marginal
density function $f_{X}$ of $X$ is continuous and strictly positive on the
interval $J=I^{\eta}$, where $I=[ c,d] $ is a compact interval and $\eta>0$.
Then it follows under the above mentioned assumptions on $\mathcal{F}$ and $%
\mathcal{K}$ that for all $c>0$ and all sequences $0<b_{n}<1$ with $%
a_{n}^{\prime\prime}\leq b_{n}\rightarrow0$, there exists a constant $%
0<C^{\prime\prime}<\infty$ such that, 
\begin{equation*}
\limsup_{n\rightarrow\infty}\sup_{a_{n}^{\prime\prime}\leq
h<b_{n}}\sup_{\varphi\in\mathcal{F}}\sup_{\mathbf{t}\in I^{m}}\frac{\sqrt{%
nh^{m}}|\hat{m}_{n, \varphi}(\mathbf{t},h)-\widehat{\mathbb{E}}\hat{m}_{n,
\varphi }(\mathbf{t},h)|}{\sqrt{|\log h|\vee\log\log n}}\leq
C^{\prime\prime},\quad\text{a.s.},
\end{equation*}
where $I^{m}=I\times\ldots\times I$ and $a_{n}^{\prime\prime}$ is either $%
a_{n}$ or $a_{n}^{\prime}$ depending on whether the class $\mathcal{F}$ is
bounded or not, i.e. whether (\ref{bounded}) or (\ref{unbounded})
holds.\medskip
\end{theo}

\noindent The following proposition follows straightforwardly from Theorems %
\ref{theo:cons} and \ref{theo:conv-rates}.

\begin{prop}
Under the assumptions of Theorems 3 and 4 on $f_{X}$ and the classes $%
\mathcal{F}$ and $\mathcal{K}$, it follows that for all sequences $0<%
\widetilde{a}_{n}\leq b_{n}<1$ satisfying $b_{n}\rightarrow0$ and $n%
\widetilde{a}_{n}/\log n\rightarrow\infty$, 
\begin{equation}
\sup_{\widetilde{a}_{n}\leq h<b_{n}}\sup_{\varphi\in\mathcal{F}}\sup_{%
\mathbf{t}\in I^{m}}|\hat{m}_{n,\varphi}(\mathbf{t},h)-m_{\varphi }(\mathbf{t%
})|\longrightarrow0,\quad a.s.,  \label{pp}
\end{equation}
where $I^{m}=I\times\ldots\times I$.
\end{prop}

\noindent It is readily seen that one can take $\widetilde{a}%
_{n}=a_{n}^{\prime}$ in the previous proposition and obtain strong uniform
consistency of Stute's estimator (\ref{def:estimator}) for general
bandwidths. However, note that by choosing $\widetilde{a}_{n}=a_{n}$, one
would only obtain almost sure convergence to a positive constant $\widetilde{%
c}>0$ in (\ref{pp}).

\section{Preliminaries for the proofs of the theorems}

Let $\Psi $ be a real valued functional defined on a class of functions $%
\mathcal{G}$ and $g$ a real valued function defined on $\mathbb{R}^d, d\geq 1$. Occasionally we shall use the notation 
\begin{equation}
\| \Psi ( G) \| _{\mathcal{G}}=\sup_{G\in \mathcal{G}}| \Psi( G) | \quad \textrm{and}\quad \| g\|_\infty =\sup_{\mathbf{x}\in \mathbb{R}^d}| g( \mathbf{x})|.  \label{TT}
\end{equation}%
In the sequel we shall need to symmetrize the functions $G_{g,h,\mathbf{t}%
}(\cdot,\cdot)$. To do this, we set 
\begin{equation*}
\bar{G}_{g,h,\mathbf{t}}(\mathbf{x},\mathbf{y}):=(m!)^{-1}\sum_{\sigma\in
I_{m}^{m}}G_{g,h,\mathbf{t}}(\mathbf{x}_{\sigma},\mathbf{y}_{\sigma}) = (m!)^{-1}\sum_{\sigma\in
I_{m}^{m}}g(\mathbf{y}_{\sigma})\widetilde{K}_{h}(\mathbf{t}-\mathbf{x}_{\sigma}),
\end{equation*}
where $\mathbf{z}_{\sigma}:=(z_{\sigma_{1}},\ldots,z_{\sigma_{m}})$.
Obviously, the expectation of $G_{g,h,\mathbf{t}}$ remains unchanged after
symmetrization, and $U_{n}^{(m)}(\bar{G}_{g,h,\mathbf{t}}(\cdot,\cdot
))=U_{n}(g,h,\mathbf{t})$, and thus the $U$--statistic process in (\ref{LU})
may be redefined using the symmetrized kernels, i.e. we consider 
\begin{equation}
u_{n}(g,h,\mathbf{t})=\sqrt{n}\{U_{n}^{(m)}(\bar{G}_{g,h,\mathbf{t}})-%
\mathbb{E}U_{n}^{(m)}(\bar{G}_{g,h,\mathbf{t}})\}.  \label{a}
\end{equation}
Moreover, the Hoeffding decomposition tells us that 
\begin{equation}
u_{n}(g,h,\mathbf{t})=\sqrt{n}\sum_{k=1}^{m}\binom{m}{k}U_{n}^{(k)}(\pi _{k}%
\bar{G}_{g,h,\mathbf{t}}(\cdot,\cdot)),  \label{hoeff:dec}
\end{equation}
where the $k$--th Hoeffding projection for the (symmetric) function $%
L:S^{m}\times S^{m}\rightarrow\mathbb{R}$ is defined for $\mathbf{x}%
_{k}=(x_{1},\ldots,x_{k})\in S^{k}$ and $\mathbf{y}_{k}=(y_{1},\ldots
,y_{k})\in S^{k}$ as 
\begin{equation*}
\pi_{k}L(\mathbf{x}_{k},\mathbf{y}_{k}):=(\delta_{(x_{1},y_{1})}-P)\times%
\ldots\times(\delta_{(x_{k},y_{k})}-P)\times P^{m-k}(L),
\end{equation*}
where $P$ is any probability measure on $(S,\mathcal{S})$. Considering $%
(X_{i},Y_{i}),i\geq1$ i.i.d--$P$ and assuming $L$ is in $L_{2}(P^{m})$, this
is an orthogonal decomposition, and $\mathbb{E} [\pi_{k}L(\mathbf{X}_{k},\mathbf{Y}_{k})|(X_{2},Y_{2}), \ldots ,(X_{k},Y_{k})]=0$,  $k\geq1$,
where we denote $\mathbf{X}_{k}$ and $\mathbf{Y}_{k}$ for $(X_{1},\ldots
,X_{k})$ and $(Y_{1},\ldots,Y_{k})$ respectively. Thus the kernels $\pi_{k}L$
are canonical for $P$ (or completely degenerate, or completely centered).
Also, $\pi_{k}$, $k\geq1$, are nested projections, i.e., $\pi_{k}\circ\pi
_{l}=\pi_{k}$ if $k\leq l$, and 
\begin{equation}
\mathbb{E}[(\pi_{k}L)^{2}(\mathbf{X}_{k},\mathbf{Y}_{k})]\leq\mathbb{E}[(L-%
\mathbb{E}L)^{2}(\mathbf{X},\mathbf{Y})]\leq\mathbb{E}L^{2}(\mathbf{X},%
\mathbf{Y}).  \label{proj}
\end{equation}
For more details consult de la Pe\~{n}a and Gin\'{e} \cite{dlPG99}.\newline

Since we assume $\mathcal{F}$ to be of VC--type with envelope function $F$,
and $\mathcal{K}$ to be of VC--type with envelope  $\kappa$, it is
readily checked (via Lemma A.1 in Einmahl and Mason \cite{EM2000}) that the
class of functions on $\mathbb{R}^{m}\times\mathbb{R}^{m}$ given by $%
\{h^{m}G_{g,h,\mathbf{t}}(\cdot,\cdot):g\in\mathcal{F},0<h<1,\mathbf{t}\in%
\mathbb{R}^{m}\}$ is of VC--type, as well as the class
\begin{equation}
\mathcal{G}=\{h^{m}\bar{G}_{g,h,\mathbf{t}}(\cdot,\cdot):g\in\mathcal{F}%
,0<h<1,\mathbf{t}\in\mathbb{R}^{m}\}, \label{gg}
\end{equation}
for which we denote the VC--type characteristics by $A_{1}$ and $v_{1}$, and
the envelope function by 
\begin{equation}
\widetilde{F}(\mathbf{y})\equiv\widetilde{F}(\mathbf{x},\mathbf{y})=\kappa
^{m}\sum_{\sigma\in I_{m}^{m}}F(\mathbf{y}_{\sigma}),\quad\mathbf{y}\in%
\mathbb{R}^{m}.  \label{envelope}
\end{equation}
(Recall $(F.ii)$ and $(F.iii)$ for terminology.) Next, for $k=1,\dots,m$
introduce the classes of functions on $\mathbb{R}^{k}\times\mathbb{R}^{k}$, 
\begin{equation}
\mathcal{G}^{(k)}=\{h^{m}\pi_{k}\bar{G}_{g,h,\mathbf{t}}(\cdot,\cdot ):g\in%
\mathcal{F},0<h<1,\mathbf{t}\in\mathbb{R}^{m}\}.  \label{G}
\end{equation}
Then an argument in Gin\'{e} and Mason \cite{GM07b} shows that each class $%
\mathcal{G}^{(k)}$ is of VC--type with characteristics $A_{1}$ and $v_{1}$
and envelope function $F_{k}\leq2^{k}\Vert\widetilde{F}\Vert_{\infty}$. (See
the completion of the proof of Theorem 1 in that paper for more details.)

\section{Proof of Theorem \protect\ref{theo:bounded} : the bounded case}

We begin with studying the first term of (\ref{hoeff:dec}), namely 
\begin{equation*}
m\sqrt{n}U_{n}^{(1)}(\pi_{1}\bar{G}_{g,h,\mathbf{t}}(\cdot,\cdot))=\frac {m}{%
\sqrt{n}}\sum_{i=1}^{n}\pi_{1}\bar{G}_{g,h,\mathbf{t}}(X_{i},Y_{i}).
\end{equation*}

\noindent {\bf {\normalsize Linear term of (\protect\ref{hoeff:dec})}}\medskip

\noindent From the definition of the Hoeffding projections and recalling that the
sample $(X_{1},Y_{1}),\ldots,(X_{n},Y_{n})$ is i.i.d., we can say for all $%
(x,y)\in\mathbb{R}^{2}$ that 
\begin{eqnarray*}
\pi_{1}\bar{G}_{g,h,\mathbf{t}}(x,y) & =&\mathbb{E}[\bar{G}_{g,h,\mathbf{t}%
}((x,X_{2},\ldots,X_{m}),(y,Y_{2},\ldots,Y_{m}))] - \mathbb{E} \bar{G}_{g,h,\mathbf{t}}(\mathbf{X},\mathbf{Y}) \\
& =&\mathbb{E}[\bar{G}_{g,h,\mathbf{t}}(\mathbf{X},\mathbf{Y}%
)|(X_{1},Y_{1})=(x,y)]-\mathbb{E}\bar{G}_{g,h,\mathbf{t}}(\mathbf{X},\mathbf{%
Y}).
\end{eqnarray*}
Introduce therefore the function on $\mathbb{R}\times\mathbb{R}$ (for
clarity we do not indicate the dependence on $m$) 
\begin{equation*}
\begin{array}{cccl}
S_{g,h,\mathbf{t}}: & \mathbb{R}\times\mathbb{R} & \longrightarrow & \mathbb{%
R} \\ 
& (x,y) & \longmapsto & mh^{m}\mathbb{E}[\bar{G}_{g,h,\mathbf{t}}(\mathbf{X},%
\mathbf{Y})|(X_{1},Y_{1})=(x,y)].%
\end{array}%
\end{equation*}
Then obviously these functions are symmetric. Using this notation we write 
\begin{equation*}
mh^{m}\pi_{1}\bar{G}_{g,h,\mathbf{t}}(x,y)=S_{g,h,\mathbf{t}}(x,y)-\mathbb{E}%
S_{g,h,\mathbf{t}}(X_{1},Y_{1}),
\end{equation*}
and hence for all $g\in\mathcal{F}$, $h\in[a_{n},b_{0}]$ and $\mathbf{t}\in%
\mathbb{R}^{m}$, the linear term of the decomposition in (\ref{hoeff:dec})
times $h^{m}$ is given by 
\begin{eqnarray*}
mh^{m}\sqrt{n}U_{n}^{(1)}(\pi_{1}\bar{G}_{g,h,\mathbf{t}}) & =&\frac{1}{\sqrt{%
n}}\sum_{i=1}^{n}\{S_{g,h,\mathbf{t}}(X_{i},Y_{i})-\mathbb{E}S_{g,h,\mathbf{t%
}}(X_{i},Y_{i})\}  \notag \\
& =:&\alpha_{n}(S_{g,h,\mathbf{t}}),  %\label{emp-proc}
\end{eqnarray*}
where this last expression is an empirical process $\alpha_{n}$ based on the
sample $(X_{1},Y_{1}),$ $\ldots,(X_{n},Y_{n})$ and indexed by the class of
functions on $\mathbb{R}\times\mathbb{R}$, 
\begin{equation*}
\mathcal{S}_{n}=\{S_{g,h,\mathbf{t}}(\cdot,\cdot):g\in\mathcal{F},a_{n}\leq
h\leq b_{0},\mathbf{t}\in\mathbb{R}^{m}\}.  \label{ss}
\end{equation*}
Clearly $\mathcal{S}_{n}\subset m\mathcal{G}^{(1)}$, and the class $m%
\mathcal{G}^{(1)}$ has envelope function $mF_{1},$ where $F_{1}$ is the
envelope function of the class $\mathcal{G}^{(1)}$ defined in (\ref{G}). From the above discussion, this class is of VC--type with the same
characteristics as $\mathcal{G}$, and therefore, after appropriate
identifications of notation, we can apply Theorem 2 of Dony, Einmahl and
Mason \cite{DEM} to conclude that 
\begin{equation}
\limsup_{n\rightarrow\infty}\sup_{a_{n}\leq h\leq b_{0}}\sup_{g\in\mathcal{F}%
}\sup_{\mathbf{t}\in\mathbb{R}^{m}}\frac{m\sqrt{nh^{m}}|U_{n}^{(1)}(\pi _{1}%
\bar{G}_{g,h,\mathbf{t}})|}{\sqrt{|\log h|\vee\log\log n}}\leq C,\quad\text{%
a.s.}  \label{conv:1st:bound}
\end{equation}
Alternatively, a straightforward modification of the proof of (\ref{lil1})
below with $a_{n}^{\prime}$ replaced by $a_{n}$ and $\gamma_{\ell}^{1/p}$ by 
$M,$ gives (\ref{conv:1st:bound}) as well.\\

\noindent {\bf {\normalsize The other terms of (\protect\ref{hoeff:dec})}}\medskip

\noindent Our aim now is to show that all the other terms of the Hoeffding
decomposition are almost surely bounded or more precisely that for each $%
k=2,\ldots,m$, 
\begin{equation}
\sup_{a_{n}\leq h\leq b_{0}}\sup_{g\in\mathcal{F}}\sup_{\mathbf{t}\in\mathbb{%
R}^{m}}\frac{\binom{m}{k}\sqrt{nh^{m}}|U_{n}^{(k)}(\pi_{k}\bar {G}_{g,h,%
\mathbf{t}})|}{\sqrt{|\log h|\vee\log\log n}}=O(1),\quad\text{a.s.}
\label{other:proj}
\end{equation}
Since $na_{n}^{m}=c^{m}\log n$, this will be accomplished if we can prove
that for each $k=2,\ldots,m$, 
\begin{equation}
\sup_{a_{n}\leq h\leq b_{0}}\sup_{g\in\mathcal{F}}\sup_{\mathbf{t}\in\mathbb{%
R}^{m}}\frac{\sqrt{nh^{m}}|U_{n}^{(k)}(\pi_{k}\bar{G}_{g,h,\mathbf{t}})|}{%
\sqrt{(|\log h|\vee\log\log n)^{k}}}=O\left( \frac {1}{\sqrt{a_{n}^{m}n^{k-1}%
}}\right) ,\quad\text{a.s.}  \label{other:projb}
\end{equation}

\noindent To obtain uniform in bandwidth convergence rates, we shall need a
blocking argument and a decomposition of the interval $[a_{n},b_{0}]$ into
smaller intervals. To do this, set $n_{\ell}=2^{\ell},\ell\geq0$ and
consider the intervals $\mathcal{H}_{\ell,j}:=[h_{\ell,j-1},h_{\ell,j}]$,
where the boundaries are given by $h_{\ell,j}^{m}:=2^{j}a_{n_{\ell}}^{m}$.
Setting $L(\ell)=\max\{j:h_{\ell,j}\leq2b_{0}\}$, observe that 
\begin{equation}
\lbrack a_{n_{\ell}},b_{0}]\subseteq\ \bigcup_{\ell=1}^{L(\ell)}\mathcal{H}%
_{\ell,j}\quad\text{and}\quad L(\ell)\sim\log\left( \frac{n_{\ell }b_{0}}{%
c\log n_{\ell}}\right) /\log2,  \label{LL}
\end{equation}
implying in particular that $L(\ell)\leq2\log n_{\ell}$. (This fact will be
used repeatedly to finish some important steps of the proofs.) Next, for $%
1\leq j\leq L(\ell)$, consider the class of functions on $\mathbb{R}%
^{m}\times\mathbb{R}^{m}$, 
\begin{equation*}
\mathcal{G}_{\ell,j}:=\{h^{m}\bar{G}_{g,h,\mathbf{t}}(\cdot,\cdot ):g\in%
\mathcal{F},h\in\mathcal{H}_{\ell,j},\mathbf{t}\in\mathbb{R}^{m}\},
\end{equation*}
as well as the class on $\mathbb{R}^{k}\times\mathbb{R}^{k}$, 
\begin{equation*}
\mathcal{G}_{\ell,j}^{(k)}:=\left\{ \frac{h^{m}\pi_{k}\bar{G}_{g,h,\mathbf{t}%
}(\cdot,\cdot)}{M_{k}}:g\in\mathcal{F},h\in\mathcal{H}_{\ell,j},\mathbf{t}\in%
\mathbb{R}^{m}\right\} ,
\end{equation*}
where $M_{k}=2^{k}\kappa^{m}M$. Clearly, each class $\mathcal{G}_{\ell,j}$
is of VC--type with the same characteristics and envelope function as $%
\mathcal{G}$, and $\mathcal{G}_{\ell,j}^{(k)}$ is of VC--type with the same
characteristics as $\mathcal{G}^{(k)}$ (and thus as $\mathcal{G}$) with
envelope function $M_{k}^{-1}F_{k}$, where $F_{k}$ is the envelope function
of $\mathcal{G}^{(k)}$. Notice that from (\ref{bounded}), 
\begin{equation*}
M_{k}\geq\sup_{\mathbf{x},\mathbf{y}\in\mathbb{R}^{k}}\{|\pi_{k}\bar{G}_{g,h,\mathbf{t}%
}(\mathbf{x},\mathbf{y})|:g\in\mathcal{F},0<h<1,\mathbf{t}\in\mathbb{R}%
^{m}\},
\end{equation*}
and hence each function in $\mathcal{G}_{\ell,j}^{(k)}$ is bounded by 1.
Define now for $n_{\ell-1}<n\leq n_{\ell}$, $\ell=1,2,\dots$ 
\begin{equation}
\mathcal{U}_{n}(j,k,\ell)=n_{\ell}^{-k/2}\sup_{H\in\mathcal{G}%
_{\ell,j}^{(k)}}\Big|\sum_{\mathbf{i}\in I_{n}^{k}}H(\mathbf{X}_{\mathbf{i}},%
\mathbf{Y}_{\mathbf{i}})\Big|.  \label{S}
\end{equation}

\noindent From Theorem 4 of Gin\'{e} and Mason \cite{GM07b} (see Theorem \ref%
{Brown} in the Appendix), we get for $c=1/2,$ $r=2$ and all $x>0$ that for
any $\ell\geq1$, 
\begin{equation}
\mathbb{P}\Big\{\max_{n_{\ell-1}<n\leq n_{\ell}}\mathcal{U}_{n}(j,k,\ell )>x%
\Big\}\leq\frac{2}{x}\mathbb{P}\left\{ \mathcal{U}_{n_{\ell}}(j,k,\ell)>x/2%
\right\} ^{1/2}\mathbb{E}[\mathcal{U}_{n_{\ell}}^{2}(j,k,\ell)]^{1/2}.
\label{ineq}
\end{equation}
We shall apply an exponential inequality and a moment bound for $U$%
--statistics due to respectively de la Pe\~{n}a and Gin\'{e} \cite{dlPG99},
and Gin\'{e} and Mason \cite{GM07b}, on the class $\mathcal{G}%
_{\ell,j}^{(k)} $ to bound (\ref{ineq}). In order to use these results we
must first derive some bounds. Firstly, it is readily checked that 
\begin{equation}
\mathcal{U}_{n}(j,k,\ell)\leq n_{\ell}^{k/2}\Vert
U_{n}^{(k)}(\pi_{k}G)\Vert_{\mathcal{G}_{\ell,j}^{(k)}},  \label{**}
\end{equation}
for all $n_{\ell-1}<n\leq n_{\ell}$. (Recall the notation (\ref{TT}).)
Secondly, notice that in $(K.i)$, $K$ is assumed to be bounded by $\kappa$ and has
support in $[-1/2,1/2]$, such that by assumption (\ref{bounded}) and $%
M_{k}=2^{k}\kappa^{m}M$, for $H\in\mathcal{G}_{\ell,j}^{(k)}$ we have by (%
\ref{proj}) 
\begin{eqnarray*}
\mathbb{E}H^{2}(\mathbf{X},\mathbf{Y}) & \leq &M_{k}^{-2}h^{2m}\mathbb{E}\bar{%
G}_{g,h,\mathbf{t}}^{2}(\mathbf{X},\mathbf{Y}) \\
& =& M_{k}^{-2}\mathbb{E}\Big[g^{2}(\mathbf{Y})\widetilde{K}^{2}\Big(\frac {%
\mathbf{t}-\mathbf{X}}{h}\Big)\Big] \\
& \leq & h^{m}4^{-k}\Vert f_{X}\Vert_{\infty}^{m}.
\end{eqnarray*}
For $D_{m}=4^{-k}\Vert f_{X}\Vert_{\infty}^{m}$, this gives us that 
\begin{equation}
\sup_{H\in\mathcal{G}_{\ell,j}^{(k)}}\mathbb{E}H^{2}(\mathbf{X},\mathbf{Y}%
)\leq D_{m}h_{\ell,j}^{m}=:\sigma_{\ell,j}^{2}.  \label{sigma}
\end{equation}
Since $\pi_{k}\pi_{k}L=\pi_{k}L$ for all $k\geq1$, we can now apply Theorem %
\ref{gm2007b} to the class $\mathcal{G}_{\ell,j}^{(k)}$ with $\sigma_{\ell
,j}^{2}$ as in (\ref{sigma}), and obtain easily that for some constant $%
A_{k} $, 
\begin{equation}
\mathbb{E}\mathcal{U}_{n_{\ell}}^{2}(j,k,\ell)\leq n_{\ell}^{k}\mathbb{E}%
\Vert U_{n_{\ell}}^{(k)}(\pi_{k}H)\Vert_{\mathcal{G}_{\ell,j}^{(k)}}^{2}\leq
2^{k}A_{k}h_{\ell,j}^{m}|\log h_{\ell,j}|^{k}.  \label{EX}
\end{equation}

\noindent To control the probability term in (\ref{ineq}), we shall apply an
exponential inequality to the same class $\mathcal{G}_{\ell,j}^{(k)}$
(recall that each $H\in\mathcal{G}_{\ell,j}^{(k)}$ is bounded by $1$).
Setting 
\begin{equation}
y^{\ast}=C_{1,k}(|\log h_{\ell,j}|\vee\log\log
n_{\ell})^{k/2}=:C_{1,k}\lambda_{j,k}(\ell),  \label{yy}
\end{equation}
where $C_{1,k}<\infty$, Theorem \ref{pena} gives us constants $%
C_{2,k},C_{3,k}$ such that for $j=1,\ldots,L(\ell)$ and any $\rho>1$, 
\begin{eqnarray}
\mathbb{P}\left\{ \mathcal{U}_{n_{\ell}}(j,k,\ell)>\rho^{k/2}y^{\ast
}\right\} & \leq &C_{2,k}\exp\left\{ -C_{3,k}\rho y^{\ast2/k}\right\}  \notag
\\
& \leq&\exp\left\{ -C_{4,k}\rho\log\log n_{\ell}\right\} .  \label{ine}
\end{eqnarray}
Then plugging the bounds (\ref{EX}) and (\ref{ine}) into (\ref{ineq}), we
get for some $C_{5,k}>0$, any $\rho\geq2$ and $\ell$ large enough, 
\begin{eqnarray}
\mathbb{P}\left\{ \max_{n_{\ell-1}<n\leq n_{\ell}}\mathcal{U}_{n}\left(
j,k,\ell\right) >2\rho^{k/2}y^{\ast}\right\} & \leq& \frac{(\log n_{\ell
})^{-\rho\frac{C_{4,k}}{2}}\sqrt{2^{k}A_{k}h_{\ell,j}^{m}|\log
h_{\ell,j}|^{k}}}{C_{1,k}\sqrt{\rho^{k}(|\log h_{\ell,j}|\vee\log\log
n_{\ell})^{k}}}  \notag \\
& \leq &\sqrt{h_{\ell,j}^{m}}(\log n_{\ell})^{-\rho C_{5,k}}.  \label{bound:GM}
\end{eqnarray}
Finally, note also that 
\begin{equation}
n_{\ell}^{k/2}\Vert U_{n}^{(k)}(\pi_{k}G)\Vert_{\mathcal{G}_{\ell,j}}\leq
C_{k}M_{k}\mathcal{U}_{n}(j,k,\ell),  \label{***}
\end{equation}
for some $C_{k}>0$. Therefore by (\ref{LL}), for each $k=2,\ldots,m$ and $%
\ell$ large enough, 
\begin{eqnarray*}
\max_{n_{\ell-1}<n\leq n_{\ell}}A_{n,k}&:= & \max_{n_{\ell-1}<n\leq
n_{\ell}}\sup_{a_{n}\leq h\leq b_{0}}\sup_{g\in\mathcal{F}}\sup_{\mathbf{t}%
\in\mathbb{R}^{m}}\frac{\sqrt{nh^{m}}|U_{n}^{(k)}(\pi_{k}\bar{G}_{g,h,%
\mathbf{t}})|}{\sqrt{(|\log h|\vee\log\log n)^{k}}} \\
& \leq& \max_{n_{\ell-1}<n\leq n_{\ell}}\max_{1\leq j\leq L(\ell)}\sup _{h\in%
\mathcal{H}_{\ell,j}}\sup_{g\in\mathcal{F}}\sup_{\mathbf{t}\in \mathbb{R}%
^{m}}\frac{\sqrt{n_{\ell}h^{m}}|U_{n}^{(k)}(\pi_{k}\bar {G}_{g,h,\mathbf{t}%
})|}{\sqrt{(|\log h|\vee\log\log n_{\ell})^{k}}} \\
& \leq &\frac{C_{k}M_{k}}{\sqrt{a_{n_{\ell}}^{m}n_{\ell}^{k-1}}}\max
_{n_{\ell-1}<n\leq n_{\ell}}\max_{1\leq j\leq L(\ell)}\frac{\mathcal{U}%
_{n}(j,k,\ell)}{\lambda_{j,k}(\ell)},
\end{eqnarray*}
where $\lambda_{j,k}(\ell)$ was defined as in (\ref{yy}). Now recall that $%
h_{\ell,j}\leq2b_{0}<2$ for $j=1,\ldots,L(\ell)$ and that $L(\ell)\leq2\log
n_{\ell}$. Then (\ref{bound:GM}) applied with $\rho\geq(2+\delta)/C_{5,k}$, $%
\delta>0$ and in combination with the above inequality and the obvious bound 
$\sqrt{a_{n}^{m}n^{k-1}}A_{n,k}\leq \sqrt{a_{n_{\ell}}^{m}n_{\ell}^{k-1}}A_{n,k}$ valid for all $n_{\ell-1}<n\leq n_{\ell}$, implies for $C_{6,k}\geq2\rho^{k/2}C_{k}M_{k}C_{1,k}$ that for $k=2,\dots,m$ 
\begin{eqnarray*}
\mathbb{P}\left\{ \max_{n_{\ell-1}<n\leq n_{\ell}}\sqrt{a_{n}^{m}n^{k-1}}%
A_{n,k}>C_{6,k}\right\} & \leq &\sum_{j=1}^{L(\ell)}\sqrt{h_{\ell,j}^{m}}(\log
n_{\ell})^{-\rho C_{5,k}} \\
& \leq &L(\ell)\sqrt{2^{m}}(\log n_{\ell})^{-\rho C_{5,k}} \\
& \leq &\sqrt{2^{m+2}}(\ell\log2)^{-(1+\delta)}.
\end{eqnarray*}
This proves via some elementary bounds and Borel--Cantelli that (\ref%
{other:projb}) holds, which obviously implies (\ref{other:proj}), and hence
completes the proof of Theorem \ref{theo:bounded}.

\section{Proof of Theorem 2 : the unbounded case}

In case (\ref{bounded}) is not satisfied, we consider bandwidths lying in
the slightly smaller interval $\mathcal{H}_{n_{\ell}}^{\prime}=[a_{n_{%
\ell}}^{\prime},b_{0}]$ that can be decomposed into the subintervals 
\begin{equation}
\mathcal{H}_{\ell,j}^{\prime}:=[h_{\ell,j-1}^{\prime},h_{\ell,j}^{\prime
}]\quad\text{with}\;\;h_{\ell,j}^{\prime m}:=2^{j}a_{n_{\ell}}^{\prime m}.
\label{H2}
\end{equation}
Note that it is straightforward to show that (\ref{LL}) remains valid if we
replace $h_{\ell,j}$ by $h_{\ell,j}^{\prime}$. In particular, we still have $%
L(\ell)\leq2\log n_{\ell}$ where $L(\ell)$ is now defined as $L(\ell
):=\max\{j:h_{\ell,j}^{\prime}\leq2b_{0}\}$. Recall that $n_{\ell}=2^{\ell
},\ell\geq0$ and set for $\ell\geq1$ 
\begin{equation}
\gamma_{\ell}=n_{\ell}/\log n_{\ell}.  \label{gamm}
\end{equation}
For an arbitrary $\varepsilon>0$ we shall decompose each function in $%
\mathcal{G}$ as 
\begin{eqnarray*}
\bar{G}_{g,h,\mathbf{t}}(\mathbf{x},\mathbf{y}) & = &\bar{G}_{g,h,\mathbf{t}}(%
\mathbf{x},\mathbf{y})\mathrm{1\!\!I}\{\widetilde{F}(\mathbf{y}%
)\leq\varepsilon\gamma_{\ell}^{1/p}\}+\bar{G}_{g,h,\mathbf{t}}(\mathbf{x},%
\mathbf{y})\mathrm{1\!\!I}\{\widetilde{F}(\mathbf{y})>\varepsilon\gamma
_{\ell}^{1/p}\} \\
& =: &\bar{G}_{g,h,\mathbf{t}}^{(\ell)}(\mathbf{x},\mathbf{y})+\widetilde {G}%
_{g,h,\mathbf{t}}^{(\ell)}(\mathbf{x},\mathbf{y})\text{,}
\end{eqnarray*}
where $\widetilde{F}(\mathbf{y})$ is the (symmetric) envelope function of
the class $\mathcal{G}$ as defined in (\ref{envelope}). Then $u_{n}(g,h,%
\mathbf{t})$ can be decomposed as well for any $n_{\ell-1}<n\leq n_{\ell}$,
since from (\ref{a}), 
\begin{eqnarray*}
u_{n}(g,h,\mathbf{t}) & =&\sqrt{{n}}\{U_{n}^{(m)}(\bar{G}_{g,h,\mathbf{t}%
}^{(\ell)})-\mathbb{E}U_{n}^{(m)}(\bar{G}_{g,h,\mathbf{t}}^{(\ell)})\} +\sqrt{n}\{U_{n}^{(m)}(\widetilde{G}_{g,h,\mathbf{t}}^{(\ell)})-\mathbb{E}U_{n}^{(m)}(\widetilde{G}_{g,h,\mathbf{t}}^{(\ell)})\}
\\
& =:&u_{n}^{(\ell)}(g,h,\mathbf{t})+\widetilde{u}_{n}^{(\ell)}(g,h,\mathbf{t}%
).
\end{eqnarray*}

\noindent The term $u_{n}^{(\ell)}(g,h,\mathbf{t})$ will be called the
truncated part and $\widetilde{u}_{n}^{(\ell)}(g,h,\mathbf{t})$ the
remainder part. To prove Theorem \ref{theo:unbounded} we shall apply the
Hoeffding decomposition to the truncated part and analyze each of the terms
separately, while the remainder part can be treated directly using simple
arguments based on standard inequalities. Note for further use that 
\begin{equation}
a_{n_{\ell}}^{\prime m}=c^{m}\gamma_{\ell}^{2/p-1},\quad\ell\geq1.
\label{gamma}
\end{equation}

\subsection{Truncated part}

%\subsubsection*{First term of the decomposition}

Note that from (\ref{hoeff:dec}) we need to consider the terms of $\sum
_{k=1}^{m}\binom{m}{k}U_{n}^{(k)}(\pi_{k}\bar{G}_{g,h,\mathbf{t}}^{(\ell)})$%
. We shall start with the linear term in this decomposition. Following the
same reasoning as in the previous section, we can show that $\pi_{1}\bar
{G}%
_{g,h,\mathbf{t}}^{(\ell)}$ is a centered conditional expectation, and that
the first term of (\ref{hoeff:dec}) can be written as an empirical process
based upon the sample $(X_{1},Y_{1}),\ldots,(X_{n},Y_{n})$ and indexed by
the class of functions 
\begin{equation*}
\mathcal{S}_{\ell}^{\prime}:=\left\{ S_{g,h,\mathbf{t}}^{(\ell)}(\cdot
,\cdot):g\in\mathcal{F},h\in\mathcal{H}_{n_{\ell}}^{\prime},\mathbf{t}\in%
\mathbb{R}^{m}\right\} ,
\end{equation*}
where $\mathcal{H}_{n_{\ell}}^{\prime}$ was defined in the beginning of this
section, and where 
\begin{equation*}
S_{g,h,\mathbf{t}}^{(\ell)}(x,y)=mh^{m}\mathbb{E}\left[ \bar{G}_{g,h,\mathbf{%
t}}^{(\ell)}(\mathbf{X},\mathbf{Y})\big|(X_{1},Y_{1})=(x,y)\right] .
\end{equation*}
To show that $\mathcal{S}_{\ell}^{\prime}$ is a VC--class, introduce the
class of functions of $(\mathbf{x},\mathbf{y})\in\mathbb{R}^{m}\times\mathbb{%
R}^{m}$, 
\begin{equation*}
\mathcal{C}=\left\{ h^{m}\bar{G}_{g,h,\mathbf{t}}(\mathbf{x},\mathbf{y})%
\mathrm{1\!\!I}\{\widetilde{F}(\mathbf{y})\leq c\}:g\in\mathcal{F},0<h<1,%
\mathbf{t}\in\mathbb{R}^{m},c>0\right\} .  %\label{gt}
\end{equation*}
Since both $\mathcal{G}$ as defined in (\ref{gg}) and the class of functions
of $\mathbf{y}\in\mathbb{R}^{m}$ given by $\mathcal{I}=\big\{\mathrm{1\!\!I}%
\{\widetilde{F}(\mathbf{y})\leq c\}:c>0\big\}$ are of VC--type (and note that $\mathcal{I}$ has a bounded envelope function), we can apply
Lemma A.1 in Einmahl and Mason \cite{EM2000} to conclude that $\mathcal{C}$  is of
VC--type as well. Therefore, so is the class of functions $m\mathcal{C}%
^{\left(  1\right)  }$ on $\mathbb{R}^{2}$, where $\mathcal{C}^{\left(
1\right)  }$ consists of the $\pi_{1}$-projections of the functions in the
class $\mathcal{C}$. Thus we see that $\mathcal{S}_{\ell}^{\prime}$ $\subset m\mathcal{C}^{\left(  1\right)  }$ and hence $\mathcal{S}_{\ell}^{\prime}$ is of VC--type with the same characteristics as $m\mathcal{C}^{\left(1\right)  }$. Now, to find an envelope function for  $\mathcal{S}_{\ell}^{\prime}$, set $\mathbf{t}_{j} :=(t_{1},\ldots,t_{j-1}, t_{j+1},\ldots,t_{m}) \in\mathbb{R}^{m-1}$, and $\mathbf{Z}_{j}(u):= (Z_{1},\ldots,Z_{j-1},u,Z_{j+1},\ldots,Z_{m})\in
\mathbb{R}^{m}$ for $u\in\mathbb{R}$ and $\mathbf{Z}\in\mathbb{R}^{m}$. We can then rewrite the function $S_{g,h,\mathbf{t}}^{(\ell)}(x,y)\in\mathcal{S}_{\ell}^{\prime}$ as
\begin{align}
S_{g,h,\mathbf{t}}^{(\ell)}(x,y)  &  =K\left(  \frac{t_{1}-x}%
{h}\right)  \mathbb{E}\Big[  g(\mathbf{Y}_{1}(y))\widetilde{K}\Big(
\frac{\mathbf{t}_{1}-\mathbf{X}^{\ast}}{h}\Big) \mathrm{1\!\!I}\{\widetilde{F} (\mathbf{Y}_{1}(y))\leq \varepsilon\gamma_{\ell}^{1/p}\}\Big] \nonumber\\
& \quad \quad +\;K\Big(  \frac{t_{2}-x}{h}\Big)  \mathbb{E}\Big[  g(\mathbf{Y}%
_{2}(y))\widetilde{K}\Big(  \frac{\mathbf{t}_{2}-\mathbf{X}^{\ast}}%
{h}\Big)  \mathrm{1\!\!I}\{\widetilde{F}(\mathbf{Y}_{2}(y))\leq \varepsilon\gamma_{\ell}^{1/p}\}\Big] \nonumber\\
&  \quad \quad+\ldots+\;K\Big(  \frac{t_{m}-x}{h}\Big)  \mathbb{E}\Big[
g(\mathbf{Y}_{m}(y))\widetilde{K}\Big(  \frac{\mathbf{t}_{m}%
-\mathbf{X}^{\ast}}{h}\Big)   \mathrm{1\!\!I}\{\widetilde{F}(\mathbf{Y}%
_{m}(y))\leq \varepsilon\gamma_{\ell}^{1/p}\}\Big]  , \label{SS}
\end{align}
where $\mathbf{X}^{\ast}=(X_{2},\ldots,X_{m})\in\mathbb{R}^{m-1}$ and where
(with abuse of notation here) the product kernel in $(K.iii)$ is now defined for
$(m-1)$--dimensional vectors, i.e. $\widetilde{K}(\mathbf{u})=\prod
_{i=1}^{m-1}K(u_{i})$, $\mathbf{u}\in\mathbb{R}^{m-1}$. Hence, we can bound
$S_{g,h,\mathbf{t}}^{(\ell)}(x,y)$ simply as
\begin{eqnarray*}
\vert S_{g,h,\mathbf{t}}^{(\ell)}(x,y)\vert &  \leq & \kappa^{m}\left\{
\mathbb{E}\left[  F(y,Y_{2},\ldots,Y_{m})\right]  +\mathbb{E}\left[
F(Y_{2},y,Y_{3},\ldots,Y_{m})\right]  \right. \\
&  &\hspace{2cm}+\ldots+\;\left.  \mathbb{E}\left[  F(Y_{2},\ldots,Y_{m}%
,y)\right]  \right\} \\
&  =:& G_{m}(x,y).
\end{eqnarray*}

\noindent We shall now apply the moment bound in Theorem \ref{em} to the
subclasses 
\begin{equation*}
\mathcal{S}_{\ell,j}^{\prime}:=\left\{ S_{g,h,\mathbf{t}}^{(\ell)}(\cdot,%
\cdot):g\in\mathcal{F},h\in\mathcal{H}_{\ell,j}^{\prime},\mathbf{t}\in%
\mathbb{R}^{m}\right\} ,\quad1\leq j\leq L(\ell),
\end{equation*}
where $\mathcal{H}_{\ell,j}^{\prime}$ was defined in (\ref{H2}). Since $%
\mathcal{S}_{\ell,j}^{\prime}\subset$ $\mathcal{S}_{\ell}^{\prime}$ for $%
j=1,\ldots,L(\ell)$, all these subclasses are of VC--type with the same
envelope function and characteristics as the class $m\mathcal{C}^{\left(
1\right) }$ (which is independent of $\ell$), verifying $(ii)$ in the
Theorem. For $(i)$, recall that although all the terms of the envelope
function $G_{m}(x,y)$ are different, their expectation is the same.
Therefore, denoting $\mathbf{Y}^{\ast}$ for $(Y_{2},\ldots,Y_{m})$ and
applying Minkowski's inequality followed by Jensen's inequality, we obtain
from assumption (\ref{unbounded}) the following upper bound for the second
moment of the envelope function.
\begin{eqnarray*}
\mathbb{E}G_{m}^{2}(X,Y) & =&\kappa^{2m}\mathbb{E}_{Y}\left\{ \mathbb{E}_{%
\mathbf{Y}^{\ast}}[F(Y,Y_{2},\ldots,Y_{m})]+\mathbb{E}_{\mathbf{Y}%
^{\ast}}[F(Y_{2},Y,Y_{3},\ldots,Y_{m})]\right. \\
& &\left. \hspace{1.5cm}+\ldots+\;\mathbb{E}_{\mathbf{Y}^{\ast}}[F(Y_{2},%
\ldots,Y_{m},Y)]\right\} ^{2} \\
& \leq &m^{2}\kappa^{2m}\mathbb{E}F^{2}(Y_{1},\ldots,Y_{m}) \\
& \leq &m^{2}\kappa^{2m}\mu_{p}^{2/p}.
\end{eqnarray*}
Note further that by symmetry of $\widetilde{F}$, 
\begin{equation*}
\mathbb{E}\bar{G}_{g,h,\mathbf{t}}^{(\ell)}(\mathbf{X},\mathbf{Y})=h^{-m}%
\mathbb{E}[g(\mathbf{Y})\widetilde{K}\Big(\frac{\mathbf{t}-\mathbf{X}}{h}\Big)\mathrm{1\!\!I}\{%
\widetilde{F}(\mathbf{Y})\leq\varepsilon\gamma_{\ell }^{1/p}\}],
\end{equation*}
such that Jensen's inequality, the change of variable $\mathbf{u}=(\mathbf{t}%
-\mathbf{x})/h$ and the assumption in (\ref{unbounded}) give the following
upper bound for the second moment of any function in $\mathcal{S}%
_{\ell}^{\prime}$ : 
\begin{eqnarray}
\mathbb{E}(S_{g,h,\mathbf{t}}^{(\ell)}(X,Y))^{2} & \leq &m^{2}\mathbb{E}\Big[%
g^{2}(\mathbf{Y})\widetilde{K}^{2}\Big(\frac{\mathbf{t}-\mathbf{X}}{h}\Big)%
\mathrm{1\!\!I}\{\widetilde{F}(\mathbf{Y})\leq\varepsilon\gamma_{\ell
}^{1/p}\}\Big]  \notag \\
&\leq &m^{2}\kappa^{2m}h^{m}\!\!\int_{\left[ -\frac{1}{2},\frac{1}{2}\right] ^{m}}\!\!\mathbb{E}[F^{2}(\mathbf{Y})\big|\mathbf{X}=\mathbf{t}-h\mathbf{u}]f_{X}(t_{1}-hu_{1})
\ldots f_{X}(t_{m}-hu_{m})d\mathbf{u}  \notag \\
& \leq &m^{2}\kappa^{2m}\mu_{p}^{2/p}\| f_{X}\|_\infty ^{m}h^{m}.
\label{dm}
\end{eqnarray}
Therefore, with $\beta\equiv m\kappa^{m}\mu_{p}^{1/p}(1\vee \| f_{X}\|_\infty ^{m})$, our previous calculations give us that 
\begin{equation*}
\mathbb{E}G_{m}^{2}(X,Y)\leq\beta^{2}\quad\text{and}\quad\sup_{S\in \mathcal{%
S}_{\ell,j}^{\prime}}\mathbb{E}S^{2}(X,Y)\leq\beta^{2}h_{\ell ,j}^{\prime m}=:\sigma_{\ell,j}^{2},
\end{equation*}
verifying condition $(iii)$ as well. Finally, recall from (\ref{envelope})
that since $\mathcal{G}$ has envelope function $\widetilde{F}(\mathbf{y})$,
it holds for all $x,y\in\mathbb{R}$ that 
\begin{equation*}
|S_{g,h,\mathbf{t}}^{(\ell)}(x,y)|\leq m\mathbb{E}[\widetilde{F}(\mathbf{Y})%
\mathrm{1\!\!I}\{\widetilde{F}(\mathbf{Y})\leq\varepsilon
\gamma_{\ell}^{1/p}\}\big|(X_{1},Y_{1})=(x,y)]\leq m\varepsilon\gamma_{\ell
}^{1/p},
\end{equation*}
such that by taking $\varepsilon>0$ small enough, Theorem \ref{em} is now
applicable, and gives us an absolute constant $A_{1}<\infty$ for which 
\begin{eqnarray}
\mathbb{E}\Vert\sum_{i=1}^{n_{\ell}}\epsilon_{i}S(X_{i},Y_{i})\Vert _{%
\mathcal{S}_{\ell,j}^{\prime}} & \leq &A_{1}\sqrt{n_{\ell}h_{\ell ,j}^{\prime
m}|\log h_{\ell,j}^{\prime}|}  \notag \\
& \leq &A_{1}\sqrt{n_{\ell}h_{\ell,j}^{\prime m}(|\log h_{\ell,j}^{\prime
}|\vee\log\log n_{\ell})}  \notag \\
&=: &A_{1}\lambda_{j}^{\prime}(\ell),  \label{ex}
\end{eqnarray}
where $\epsilon_{1},\ldots,\epsilon_{n_{\ell}}$ are independent Rademacher
variables, independent of $(X_{i},Y_{i}),$ $1\leq i\leq n_{\ell}$.
Consequently, applying the exponential inequality of Talagrand \cite{TAL94}
to the class $\mathcal{S}_{\ell,j}^{\prime}$ (see Theorem \ref{tala} in the
Appendix) with $M=m\varepsilon\gamma_{\ell}^{1/p}$, $\sigma_{\mathcal{S}%
_{\ell,j}^{\prime}}^{2}=\beta^{2}h_{\ell,j}^{\prime m}$ and the moment bound
in (\ref{ex}), we get for an absolute constant $A_{2}<\infty$ and all $t>0$
that 
\begin{equation*}
\hspace{-2.5cm}\mathbb{P}\left\{ \max_{n_{\ell-1}<n\leq n_{\ell}}\Vert \sqrt{%
n}\alpha_{n}\Vert_{\mathcal{S}_{\ell,j}^{\prime}}\geq
C_{1}(A_{1}\lambda_{j}^{\prime}(\ell)+t)\right\}
\end{equation*}%
\begin{equation}\hspace{3cm}\leq2\left[ \exp\left( -\frac{A_{2}t^{2}}{n_{\ell}\beta
^{2}h_{\ell,j}^{\prime m}}\right) +\exp\left( -\frac{A_{2}t}{m\varepsilon
\gamma_{\ell}^{1/p}}\right) \right] .  \label{TA}
\end{equation}
Towards applying this inequality with $t=\rho\lambda_{j}^{\prime}(\ell
),\rho>1$, note that it clearly follows from (\ref{gamma}) and the
definitions of $h_{\ell,j}^{\prime}$ and $\lambda_{j}^{\prime}(\ell)$ that
for all $j\geq0$, 
\begin{align*}
& \frac{\lambda_{j}^{\prime2}(\ell)}{n_{\ell}h_{\ell,j}^{\prime m}}=|\log
h_{\ell,j}^{\prime}|\vee\log\log n_{\ell}\geq\log\log n_{\ell}, \\
& \frac{\lambda_{j}^{\prime2}(\ell)}{\gamma_{\ell}^{2/p}}=2^{j}c^{m}\log
n_{\ell}(|\log h_{\ell,j}^{\prime}|\vee\log\log n_{\ell})\geq c^{m}(\log\log
n_{\ell})^{2}.
\end{align*}
Consequently, (\ref{TA}) when applied with $t=\rho\lambda_{j}^{\prime}(\ell)$
and any $\rho>1$ with $\ell$ large enough, yields for suitable constants $%
A_{2}^{\prime}$, $A_{2}^{\prime\prime}$ and $A_{3}$, the inequality 
\begin{equation*}
\hspace{-4cm}\mathbb{P}\left\{ \max_{n_{\ell-1}<n\leq n_{\ell}}\Vert\sqrt {n}%
\alpha_{n}\Vert_{\mathcal{S}_{\ell,j}^{\prime}}\geq
C_{1}(A_{1}+\rho)\lambda_{j}^{\prime}(\ell)\right\} \vspace{-5mm}
\end{equation*}%
\begin{eqnarray}
\hspace{2cm} & \leq &2\left[ \exp\left( -A_{2}^{\prime}\rho^{2}\log\log
n_{\ell}\right) +\exp\left( -A_{2}^{\prime\prime}\rho\log\log n_{\ell
}\right) \right]  \notag \\
& \leq &4(\log n_{\ell})^{-A_{3}\rho}.  \label{bound:tala}
\end{eqnarray}
Keeping in mind that $mh^{m}\sqrt{n}U_{n}^{(1)}(\pi_{1}\bar{G}_{g,h,\mathbf{t%
}}^{(\ell)})$ is an empirical process $\alpha_{n}(S_{g,h,\mathbf{t}%
}^{(\ell)})$ indexed by the class $\mathcal{S}_{\ell}^{\prime}$, and
recalling (\ref{LL}), we obtain for $\ell\geq1$ that, 
\begin{equation*}
\hspace{-2.5cm}\max_{n_{\ell-1}<n\leq n_{\ell}}A_{n,\ell}^{\prime}:=\max_{n_{\ell-1}<n\leq
n_{\ell}}\sup_{a_{n}^{\prime}\leq h\leq b_{0}}\sup_{g\in\mathcal{F}}\sup_{%
\mathbf{t}\in\mathbb{R}^{m}}\frac{m\sqrt{nh^{m}}|U_{n}^{(1)}(\pi_{1}\bar{G}%
_{g,h,\mathbf{t}}^{(\ell)})|}{\sqrt{|\log h|\vee\log\log n}}\vspace{-2mm}
\end{equation*}%
\begin{eqnarray*}
\hspace{2cm}& \leq &\max_{n_{\ell-1}<n\leq n_{\ell}}\max_{1\leq j\leq L(\ell)}\sup _{h\in%
\mathcal{H}_{\ell,j}^{\prime}}\sup_{g\in\mathcal{F}}\sup_{\mathbf{t}\in%
\mathbb{R}^{m}}\frac{2\sqrt{2}\:|\sqrt{n}\alpha_{n}(S_{g,h,\mathbf{t}%
}^{(\ell)})|}{\sqrt{n_{\ell}h_{\ell,j}^{\prime m}(|\log h_{\ell,j}^{\prime
}|\vee\log\log n_{\ell})}} \\
& \leq& \max_{n_{\ell-1}<n\leq n_{\ell}}\max_{1\leq j\leq L(\ell)}\sup _{H\in%
\mathcal{S}_{\ell,j}^{\prime}}\frac{3|\sqrt{n}\alpha_{n}(H)|}{%
\lambda_{j}^{\prime}(\ell)}.
\end{eqnarray*}
Consequently, recalling once again that $L(\ell)\leq2\log n_{\ell}$, we can
infer from (\ref{bound:tala}) that for some constant $C_{5}(\rho)\geq
3C_{1}(A_{1}+\rho)$, 
\begin{eqnarray*}
 \mathbb{P}\left\{ \max_{n_{\ell-1}<n\leq n_{\ell}}A_{n,\ell}^{\prime
}>C_{5}(\rho)\right\}  &\leq &\sum_{j=1}^{L(\ell)}\mathbb{P}\left\{ \max_{n_{\ell-1}<n\leq
n_{\ell}}\Vert\sqrt{n}\alpha_{n}\Vert_{\mathcal{S}_{\ell,j}^{%
\prime}}>C_{1}(A_{1}+\rho)\lambda_{j}^{\prime}(\ell)\right\} \\
& \leq &8(\log n_{\ell})^{1-A_{3}\rho}.
\end{eqnarray*}
The Borel--Cantelli lemma when combined with this inequality for $\rho
\geq(2+\delta)/A_{3}$, $\delta>0$ and with the choice $n_{\ell}=2^{\ell}$,
establish for some $C^{\prime}<\infty$ and with probability one, that 
\begin{equation}
\limsup_{\ell\rightarrow\infty}\max_{n_{\ell-1}<n\leq n_{\ell}}\sup
_{a_{n}^{\prime}\leq h\leq b_{0}}\sup_{g\in\mathcal{F}}\sup_{\mathbf{t}\in%
\mathbb{R}^{m}}\frac{m\sqrt{nh^{m}}|U_{n}^{(1)}(\pi_{1}\bar{G}_{g,h,\mathbf{t%
}}^{(\ell)})|}{\sqrt{|\log h|\vee\log\log n}}\leq C^{\prime},  \label{lil1}
\end{equation}
finishing the study of the first term in (\ref{hoeff:dec}). We now show that
all the other terms of (\ref{hoeff:dec}) are asymptotically bounded or go to
zero at the proper rate, which will be obtained if we can prove that for $%
k=2,\ldots,m$ and with probability one, 
\begin{equation}
\max_{n_{\ell-1}<n\leq n_{\ell}}\sup_{a_{n}^{\prime}\leq h\leq
b_{0}}\sup_{g\in\mathcal{F}}\sup_{\mathbf{t}\in\mathbb{R}^{m}}\frac{\sqrt{%
nh^{m}}|U_{n}^{(k)}(\pi_{k}\bar{G}_{g,h,\mathbf{t}}^{(\ell)})|}{\sqrt{|\log
h|\vee\log\log n}}=O(\gamma_{\ell}^{1-k/2}).  \label{other-terms}
\end{equation}
Analogously to the bounded case, we start by defining the classes of
functions on $\mathbb{R}^{m}\times\mathbb{R}^{m}$ and $\mathbb{R}^{k}\times 
\mathbb{R}^{k}$, 
\begin{align*}
& \mathcal{G}_{\ell,j}^{\prime}:=\left\{ h^{m}\bar{G}_{g,h,\mathbf{t}%
}^{(\ell)}(\cdot,\cdot):g\in\mathcal{F},h\in\mathcal{H}_{\ell,j}^{\prime },%
\mathbf{t}\in\mathbb{R}^{m}\right\} , \\
& \mathcal{G}_{\ell,j}^{\prime(k)}:=\left\{ h^{m}(\pi_{k}\bar{G}_{g,h,%
\mathbf{t}}^{(\ell)})(\cdot,\cdot)/(2^{k}\varepsilon\gamma_{\ell}^{1/p}):g\in%
\mathcal{F},h\in\mathcal{H}_{\ell,j}^{\prime},\mathbf{t}\in\mathbb{R}%
^{m}\right\} .
\end{align*}
Then it is easily verified that these classes are of VC--type with
characteristics that are independent of $\ell$, and with envelope functions $%
\widetilde{F}$ and $(2^{k}\varepsilon\gamma_{\ell}^{1/p})^{-1}F_{k}$
respectively. The function $\widetilde{F}$ is defined as in (\ref{envelope})
and $F_{k}$ is determined just as in the proof of Theorem 1 in Gin\'{e} and
Mason \cite{GM07b}. Note that, in the same spirit as (\ref{S}) and (\ref{**}%
), by setting 
\begin{equation*}
\mathcal{U}_{n}^{\prime}(j,k,\ell):=\sup_{H\in\mathcal{G}_{\ell,j}^{\prime
(k)}}\Big|\frac{1}{n_{\ell}^{k/2}}\sum_{\mathbf{i}\in I_{n}^{k}}H(\mathbf{X}%
_{\mathbf{i}},\mathbf{Y}_{\mathbf{i}})\Big|,\quad n_{\ell-1}<n\leq n_{\ell},
\end{equation*}
we have for all $k=2,\ldots,m$ and $n_{\ell-1}<n\leq n_{\ell}$, 
\begin{equation*}
\mathcal{U}_{n}^{\prime}(j,k,\ell)\leq n_{\ell}^{k/2}\Vert U_{n}^{(k)}(\pi
_{k}G)\Vert_{\mathcal{G}_{\ell,j}^{\prime(k)}}.
\end{equation*}
Consequently, applying Theorem \ref{Brown} with $c=1/2$ and $r=2$, gives us
precisely (\ref{ineq}) with $\mathcal{U}_{n}(j,k,\ell)$ and $\mathcal{U}%
_{n_{\ell}}(j,k,\ell)$ replaced by $\mathcal{U}_{n}^{\prime}(j,k,\ell)$ and $%
\mathcal{U}_{n_{\ell}}^{\prime}(j,k,\ell)$ respectively. Therefore the same
methodology as in the bounded case will be applied. Note also that, as held
for all the functions in $\mathcal{G}_{\ell,j}^{(k)}$, the functions in $%
\mathcal{G}_{\ell,j}^{\prime(k)}$ are bounded by 1, and have second moments
that can be bounded by $h^{m}D_{m}$ for a suitable\ $D_{m}$ by arguing as in
(\ref{dm}) and (\ref{sigma}). Consequently, the expression in (\ref{sigma})
is satisfied for functions in $\mathcal{G}_{\ell,j}^{\prime(k)}$ as well, i.e. 
\begin{equation*}
\sup_{H\in\mathcal{G}_{\ell,j}^{\prime(k)}}\mathbb{E}H^{2}(\mathbf{X},%
\mathbf{Y})\leq D_{m}h_{\ell,j}^{\prime m}=:\sigma_{\ell,j}^{\prime2}.
\end{equation*}
Hence, all the conditions for Theorems \ref{gm2007b} and \ref{pena} are
satisfied, so that after some obvious identifications and modifications, the
second part of the proof of Theorem \ref{theo:bounded} (and (\ref{bound:GM})
in particular) gives us for all $j=1,\ldots,L(\ell)$ and any $\rho>2$, 
\begin{equation}
\mathbb{P}\left\{ \max_{n_{\ell-1}<n\leq n_{\ell}}\mathcal{U}_{n}^{\prime
}(j,k,\ell)>2\rho^{k/2}y^{\prime}{}^{\ast}\right\} \leq\sqrt{h_{\ell
,j}^{\prime m}}(\log n_{\ell})^{-C_{7,k}\rho},  \label{bound:snjkl}
\end{equation}
with $y^{\prime}{}^{\ast}=C_{1,k}^{\prime}\lambda_{j,k}^{\prime}(\ell)$, and
where $\lambda_{j,k}^{\prime}(\ell)$ is defined as in (\ref{yy}) with $%
h_{\ell,j}$ replaced by $h_{\ell,j}^{\prime}$. Now, to finish the proof of (%
\ref{other-terms}), note that similarly to (\ref{***}), for some $C_{k}>0$, 
\begin{equation*}
n_{\ell}^{k/2}\Vert U_{n}^{(k)}(\pi_{k}G)\Vert_{\mathcal{G}%
_{\ell,j}^{\prime}}\leq2^{k}C_{k}\varepsilon\gamma_{\ell}^{1/p}\mathcal{U}%
_{n}^{\prime}(j,k,\ell)\text{.}
\end{equation*}
This gives that 
\begin{eqnarray*}
\max_{n_{\ell-1}<n\leq n_{\ell}}A_{n,\ell,k}^{\prime} & :=&\max_{n_{\ell
-1}<n\leq n_{\ell}}\sup_{a_{n}^{\prime}\leq h\leq b_{0}}\sup_{g\in\mathcal{F}%
}\sup_{\mathbf{t}\in\mathbb{R}^{m}}\frac{\sqrt{nh^{m}}|U_{n}^{(k)}(\pi_{k}%
\bar{G}_{g,h,\mathbf{t}}^{(\ell)})|}{\sqrt{(|\log h|\vee\log\log n)^{k}}} \\
& \leq& \frac{2^{k}c_{k}\varepsilon\gamma_{\ell}^{1/p}}{\sqrt{%
a_{n_{\ell}}^{\prime m}n_{\ell}^{k-1}}}\max_{n_{\ell-1}<n\leq
n_{\ell}}\max_{1\leq j\leq L(\ell)}\frac{\mathcal{U}_{n}^{\prime}(j,k,\ell)}{%
\lambda_{j,k}^{\prime}(\ell)}.
\end{eqnarray*}

\noindent From (\ref{gamma}) we see now that $\gamma_{\ell}^{2/p}/a_{n_{\ell}}^{\prime
m}n_{\ell}^{k-1}=c^{-m}n_{\ell}^{2-k}/\log n_{\ell}$. Therefore by choosing $%
C_{8,k}>2^{k+1}c^{-m/2}\varepsilon c_{k}C_{1,k}^{\prime}((2+\delta
)/C_{7,k})^{k/2}$ and noting that $h_{\ell,j}^{\prime}<2$ for all $%
j=1,\ldots,L(\ell)$, we can infer from (\ref{bound:snjkl}) that 
\begin{equation*}
\mathbb{P}\Big\{ \max_{n_{\ell-1}<n\leq n_{\ell}}\sqrt{\frac{\log n}{n^{2-k}%
}}A_{n,\ell,k}^{\prime}>C_{8,k}\Big\} \leq\sqrt{2^{m+1}}(\log
n_{\ell})^{-(1+\delta)}.
\end{equation*}
This implies immediately via Borel--Cantelli that for all $k=2,\ldots,m$ and 
$\ell\geq1$, 
\begin{equation*}
\max_{n_{\ell-1}<n\leq n_{\ell}}\sup_{a_{n}^{\prime}\leq h\leq
b_{0}}\sup_{g\in\mathcal{F}}\sup_{\mathbf{t}\in\mathbb{R}^{m}}\frac{\sqrt{%
nh^{m}}|U_{n}^{(k)}(\pi_{k}\bar{G}_{g,h,\mathbf{t}}^{(\ell)})|}{\sqrt{(|\log
h|\vee\log\log n)^{k}}}=O\left( \sqrt{\frac{n_{\ell}^{2-k}}{\log n_{\ell}}}%
\right) , \quad\textrm{a.s.},
\end{equation*}
which obviously implies (\ref{other-terms}). Finally, recalling the
Hoeffding decomposition (\ref{hoeff:dec}), this implies together with (\ref%
{lil1}) that with probability one, 
\begin{equation}
\limsup_{\ell\rightarrow\infty}\max_{n_{\ell-1}<n\leq n_{\ell}}\sup
_{a_{n}^{\prime}\leq h\leq b_{0}}\sup_{g\in\mathcal{F}}\sup_{\mathbf{t}\in%
\mathbb{R}^{m}}\frac{\sqrt{nh^{m}}|U_{n}^{(m)}(\bar{G}_{g,h,\mathbf{t}%
}^{(\ell)})-\mathbb{E}U_{n}^{(m)}(\bar{G}_{g,h,\mathbf{t}}^{(\ell)})|}{\sqrt{%
|\log h|\vee\log\log n}}\leq C^{^{\prime\prime}}.  \label{result:trunc}
\end{equation}

\subsection{Remainder part}

Consider now the remainder process $\widetilde{u}_{n}^{(\ell)}(g,h,\mathbf{t}%
)$ based on the unbounded (symmetric) $U$--kernel given by 
\begin{equation*}
\widetilde{G}_{g,h,\mathbf{t}}^{(\ell)}(\mathbf{x},\mathbf{y}):=\bar {G}%
_{g,h,\mathbf{t}}(\mathbf{x},\mathbf{y})\mathrm{1\!\!I}\{\widetilde {F}(%
\mathbf{y})>\varepsilon\gamma_{\ell}^{1/p}\},
\end{equation*}
where we defined $\gamma_\ell$ as in (\ref{gamm}). We shall show that this $U$--process is
asymptotically negligible at the rate given in Theorem \ref{theo:unbounded}.
More precisely, we shall prove that as $\ell\rightarrow\infty$, 
\begin{equation}
\max_{n_{\ell-1}<n\leq n_{\ell}}\sup_{a_{n}^{\prime}\leq h\leq
b_{0}}\sup_{g\in\mathcal{F}}\sup_{\mathbf{t}\in\mathbb{R}^{m}}\frac{\sqrt{%
nh^{m}}|U_{n}^{(m)}(\widetilde{G}_{g,h,\mathbf{t}}^{(\ell)})-\mathbb{E}%
U_{n}^{(m)}(\widetilde{G}_{g,h,\mathbf{t}}^{(\ell)})|}{\sqrt{|\log
h|\vee\log\log n}}=o(1),  \quad\textrm{a.s.}\label{toprove}
\end{equation}
Recall that for all $g\in\mathcal{F},$ $h\in\lbrack a_{n}^{\prime},b_{0}]$
and $\mathbf{t},\mathbf{x}\in\mathbb{R}^{m}$, $\widetilde{F}(\mathbf{y})\geq h^{m}|\bar{G}_{g,h,\mathbf{t}}(\mathbf{x},\mathbf{y})|$, so from the symmetry of $\widetilde{F}$, it holds that 
\begin{equation*}
|U_{n}^{(m)}(\widetilde{G}_{g,h,\mathbf{t}}^{(\ell)})|\leq
h^{-m}U_{n}^{(m)}\left( \widetilde{F}\cdot\mathrm{1\!\!I}\{\widetilde{F}%
>\varepsilon \gamma_{\ell}^{1/p}\}\right) ,  %\label{ustat}
\end{equation*}
where $U_{n}^{(m)}(\widetilde{F}\cdot\mathrm{1\!\!I}\{\widetilde {F}%
>\varepsilon\gamma_{\ell}^{1/p}\})$ is a $U$--statistic based on the
positive and symmetric kernel $\mathbf{y}\rightarrow\widetilde{F}(\mathbf{y})%
\mathrm{1\!\!I}\{\widetilde{F}(\mathbf{y})>\varepsilon\gamma _{\ell}^{1/p}\}$%
. Recalling that $a_{n}^{\prime m}=c^{m}(\log n/n)^{1-2/p}$, we obtain easily that for all $g\in\mathcal{F},$ $h\in\lbrack a_{n}^{\prime
},b_{0}],$ $\mathbf{t}\in\mathbb{R}^{m}$ and some $C>0$,

\begin{eqnarray*}
\max_{n_{\ell-1}<n\leq n_{\ell}}\frac{\sqrt{nh^{m}}|U_{n}^{(m)}(\widetilde {G%
}_{g,h,\mathbf{t}}^{(\ell)})|}{\sqrt{|\log h|\vee\log\log n}} & \leq& \frac{%
\sqrt{n_{\ell}}\:U_{n_{\ell}}^{(m)}\big(\widetilde{F}\cdot\mathrm{1\!\!I}\{%
\widetilde{F}>\varepsilon\gamma_{\ell}^{1/p}\}\big)}{\sqrt{%
a_{n_{\ell}}^{\prime m}(|\log a_{n_{\ell}}^{\prime}|\vee\log\log n_{\ell})}}
\\
& \leq &C\gamma_{\ell}^{1-1/p}U_{n_{\ell}}^{(m)}\big(\widetilde{F}\cdot%
\mathrm{1\!\!I}\{\widetilde{F}>\varepsilon\gamma_{\ell}^{1/p}\}\big).
\end{eqnarray*}
Arguing in the same way, since a $U$--statistic is an unbiased estimator of
its kernel, we get  that uniformly in $g\in\mathcal{F},$ $h\in\lbrack
a_{n}^{\prime},b_{0}]$ and $\mathbf{t}\in\mathbb{R}^{m}$, 
\begin{align}
\max_{n_{\ell-1}<n\leq n_{\ell}}\frac{\sqrt{nh^{m}}|\mathbb{E}U_{n}^{(m)}(%
\widetilde{G}_{g,h,\mathbf{t}}^{(\ell)})|}{\sqrt{|\log h|\vee\log\log n}} &
\leq C\gamma_{\ell}^{1-1/p}\mathbb{E}U_{n_{\ell}}^{(m)}\big(\widetilde {F}%
\cdot\mathrm{1\!\!I}\{\widetilde{F}>\varepsilon\gamma_{\ell}^{1/p}\}\big) 
\notag \\
& \leq C^{\prime}\mathbb{E}[\widetilde{F}^{p}(\mathbf{Y})\mathrm{1\!\!I}\{%
\widetilde{F}(\mathbf{Y})>\varepsilon\gamma_{\ell}^{1/p}\}].  \label{bb}
\end{align}

\noindent From (\ref{bb}) we see that as $\ell\rightarrow\infty$, 
\begin{equation}
\max_{n_{\ell-1}<n\leq n_{\ell}}\sup_{a_{n}^{\prime}\leq h\leq
b_{0}}\sup_{g\in\mathcal{F}}\sup_{\mathbf{t}\in\mathbb{R}^{m}}\frac{\sqrt{%
nh^{m}}|\mathbb{E}U_{n}^{(m)}(\widetilde{G}_{g,h,\mathbf{t}}^{(\ell)})|}{%
\sqrt{|\log h|\vee\log\log n}}=o(1)\text{.}  \label{aa}
\end{equation}
Thus to finish the proof of (\ref{toprove}) it suffices to show that 
\begin{equation}
U_{n_{\ell}}^{(m)}(\widetilde{F}\cdot\mathrm{1\!\!I}\{\widetilde {F}%
>\varepsilon\gamma_{\ell}^{1/p}\})=o(\gamma_{\ell}^{1/p-1}),\quad \text{a.s.}
\label{ddd}
\end{equation}
First note that from Chebyshev's inequality and a well--known inequality for
the variance of a $U$--statistic (see Theorem 5.2 of Hoeffding \cite{HOEF48}%
) we get for any $\delta>0$, 
\begin{equation*}
\mathbb{P}\left\{ \big\vert U_{n_{\ell}}^{(m)}\big(\widetilde{F}\cdot%
\mathrm{1\!\!I}\{\widetilde{F}>\varepsilon\gamma_{\ell}^{1/p}\}\big)-\mathbb{%
E}U_{n_{\ell}}^{(m)}\big(\widetilde{F}\cdot\mathrm{1\!\!I}\{\widetilde{F}%
>\varepsilon\gamma_{\ell}^{1/p}\}\big)\big\vert >\delta
\gamma_{\ell}^{-(1-1/p)}\right\}
\end{equation*}
\vspace{-5mm}
\begin{eqnarray}
\hspace{1.5cm}& \leq& \delta^{-2}\gamma_{\ell}^{2-2/p}\mathrm{Var}\left( U_{n_{\ell}}^{(m)}%
\big(\widetilde{F}\cdot\mathrm{1\!\!I}\{\widetilde{F}>\varepsilon
\gamma_{\ell}^{1/p}\}\big)\right)  \notag \\
& \leq & m\delta^{-2}\frac{n_{\ell}^{1-2/p}}{(\log n_{\ell})^{2-2/p}}\mathbb{E}%
[\widetilde{F}^{2}(\mathbf{Y})\mathrm{1\!\!I}\{\widetilde {F}(\mathbf{Y}%
)>\varepsilon\gamma_{\ell}^{1/p}\}].  \label{cheby}
\end{eqnarray}
Next, in order to establish the finite convergence of the series of the
above probabilities, we split the indicator function $\mathrm{1\!\!I}\{%
\widetilde {F}(\mathbf{Y})>\varepsilon\gamma_{\ell}^{1/p}\}$ into two
distinct parts determined by whether  $\widetilde{F}(\mathbf{Y})>n_{\ell}^{1/p}$ or $\varepsilon\gamma_{\ell}^{1/p}<\widetilde{F}(\mathbf{Y})\leq n_{\ell}^{1/p}$, and consider the corresponding second moments in (\ref{cheby}) separately. In the second case, note that from (\ref{unbounded}) and (\ref{envelope}), $ \mathbb{E}\widetilde{F}^{p}(\mathbf{Y})\leq\mu_{p}\kappa^{pm}(m!)^{p}$, and
observe that since $p>2$ and $n_{\ell}=2^{\ell}$, 
\begin{equation*}
\hspace{-1cm}\sum_{\ell=1}^{\infty}\frac{n_{\ell}^{1-2/p}}{(\log n_{\ell
})^{2-2/p}}\mathbb{E}[\widetilde{F}^{2}(\mathbf{Y})\mathrm{1\!\!I}\{%
\widetilde{F}(\mathbf{Y})>n_{\ell}^{1/p}\}] \leq\mathbb{E}[\widetilde{F}^{p}(\mathbf{Y})] \sum_{\ell=1}^{\infty}(\log n_{\ell})^{-(2-2/p)}<\infty.
\end{equation*}
To handle the first case, we shall need the following fact from Einmahl and
Mason \cite{EM2000}.

\begin{fact}
Let $(c_{n})_{n\geq1}$ be a sequence of positive constants such that $%
c_{n}/n^{1/s}\nearrow\infty$ for $s>0,$ and let $Z$ be a random variable
satisfying $\sum_{n=1}^{\infty}\mathbb{P}\{|Z|>c_{n}\}<\infty$. Then we have
for any $q>s$, 
\begin{equation*}
\sum_{k=1}^{\infty}2^{k}\mathbb{E}[|Z|^{q}\mathrm{1\!\!I}\{|Z|\leq
c_{2^{k}}\}]/(c_{2^{k}})^{q}<\infty.
\end{equation*}
\end{fact}

\noindent Setting $c_{n}=n^{1/p}$ into Fact 1, we conclude from this
inequality that for $p<s<r\leq2p$, 
\begin{equation*}
\hspace{-2.5cm}\sum_{\ell=1}^{\infty}\frac{n_{\ell}^{1-2/p}}{(\log n_{\ell
})^{2-2/p}}\mathbb{E}[\widetilde{F}^{2}(\mathbf{Y})\mathrm{1\!\!I}%
\{\varepsilon\gamma_{\ell}^{1/p}<\widetilde{F}(\mathbf{Y})\leq
n_{\ell}^{1/p}\}]
\end{equation*}%
\begin{equation*}
\hspace{2cm}\leq\;\sum_{\ell=1}^{\infty}\frac{\varepsilon^{r-2}}{(\log
n_{\ell})^{2-r/p}}\frac{n_{\ell}\mathbb{E}[\widetilde{F}^{r}(\mathbf{Y})%
\mathrm{1\!\!I}\{\widetilde{F}(\mathbf{Y})\leq n_{\ell}^{1/p}\}]}{n_{\ell
}^{r/p}}<\infty.  %\label{emf1}
\end{equation*}

\noindent Finally, note that the bound leading to (\ref{bb}) implies that 
\begin{equation*}
\gamma_{\ell}^{1-1/p}\mathbb{E}U_{n_{\ell}}^{(m)}\big( \widetilde{F}\cdot%
\mathrm{1\!\!I}\{\widetilde{F}>\varepsilon\gamma_{\ell}^{1/p}\}\big) =o(1)%
\text{.}
\end{equation*}
Consequently, the above results together with (\ref{cheby}) imply via
Borel-Cantelli and the arbitrary choice of $\delta>0$ that (\ref{ddd})
holds, which when combined with (\ref{aa}) and (\ref{bb}) completes the
proof of (\ref{toprove}). This also finishes the proof of Theorem \ref%
{theo:unbounded} since we have already established the result in (\ref%
{result:trunc}).

\section{Proof of Theorem \protect\ref{theo:cons} : uniform consistency of $%
\hat
{m}_{n}(\mathbf{t},h)$ to $m_{\protect\varphi}(\mathbf{t})$}

Theorem 3 is essentially a consequence of Theorem \ref{stein} in the
Appendix. Recall that a $U$--statistic with $U$--kernel $H$ is an unbiased
estimator of $\mathbb{E}H$. Writing $d\mathbf{x}$ and $d\mathbf{y}$ for $%
dx_{1}dx_{2}\ldots dx_{m}$ and $dy_{1}dy_{2}\ldots dy_{m}$ respectively, we
see that 
\begin{align}
\mathbb{E}U_{n}(1,h,\mathbf{t}) & =\int\widetilde{K}_{h}(\mathbf{t}-\mathbf{x%
})f(x_{1},y_{1})\cdots f(x_{m},y_{m})d\mathbf{x}d\mathbf{y}   =\widetilde{f} \ast \widetilde{K}_{h}(\mathbf{t}),  \notag
\end{align}
where the function $\widetilde{f}$ $:\mathbb{R}^{m}\rightarrow\mathbb{R}$ is
defined in (\ref{ff}). Since we assume $f_{X}$ to be continuous on $%
J=I^{\eta }$, the function $\widetilde{f}$ is continuous on $%
J^{m}=J\times\ldots\times J$. Therefore we can infer from Theorem \ref{stein}
that 
\begin{equation}
\sup_{0<h<b_{n}}\sup_{\mathbf{t}\in I^{m}}|\mathbb{E}U_{n}(1,h,\mathbf{t})-%
\widetilde{f}(\mathbf{t})|\longrightarrow0,  \label{stein1}
\end{equation}
for all sequences of positive constants $b_{n}\rightarrow0$, and where $%
I^{m}=I\times\ldots\times I$. In the same way, notice that 
\begin{eqnarray*}
\mathbb{E}U_{n}(\varphi,h,\mathbf{t}) & =&\int\varphi(\mathbf{y})\widetilde{K}%
_{h}(\mathbf{t}-\mathbf{x})f(x_{1},y_{1})\cdots f(x_{m},y_{m})d\mathbf{x}d%
\mathbf{y} \\
& =&\left\{ \mathbb{E}\left[ \varphi(\mathbf{Y})|\mathbf{X}=\cdot\right] 
\text{ }\widetilde{f}(\cdot)\right\} \ast\widetilde{K}_{h}(\mathbf{t}).
\end{eqnarray*}
Hence, Theorem \ref{stein} applied to the class of functions $\mathcal{M}$
as defined in (\ref{equi}) gives that 
\begin{equation}
\sup_{0<h<b_{n}}\sup_{\varphi\in\mathcal{F}}\sup_{\mathbf{t}\in I^{m}}|%
\mathbb{E}U_{n}(\varphi,h,\mathbf{t})-m_{\varphi}(\mathbf{t})\widetilde {f}(%
\mathbf{t})|\longrightarrow0.  \label{stein2}
\end{equation}

\noindent Keeping in mind the definition of $\widehat{\mathbb{E}}\hat
{m}%
_{n,\varphi}(\mathbf{t},h)$ in (\ref{def:cent}), it is clear that since $%
f_{X}$ is bounded away from zero on $J$, (\ref{stein1}) and (\ref{stein2})
imply that 
\begin{equation*}
\sup_{0< h<b_{n}}\sup_{\varphi\in\mathcal{F}}\sup_{\mathbf{t}\in I^{m}}|%
\widehat{\mathbb{E}}\hat{m}_{n,\varphi}(\mathbf{t},h)-m_{\varphi}(\mathbf{t}%
)|=o(1),
\end{equation*}
finishing the proof of Theorem \ref{theo:cons}.

\section{Proof of Theorem \protect\ref{theo:conv-rates} : convergence rates
of the conditional $U$--statistic $\hat{m}_{n,\protect\varphi}(\mathbf{t},h)$%
}

Observe that 
\begin{eqnarray*}
|\hat{m}_{n,\varphi }(\mathbf{t},h)-\widehat{\mathbb{E}}\hat{m}_{n,\varphi }(%
\mathbf{t},h)|& =&\left\vert \frac{U_{n}(\varphi ,h,\mathbf{t})}{U_{n}(1,h,%
\mathbf{t})}-\frac{\mathbb{E}U_{n}(\varphi ,h,\mathbf{t})}{\mathbb{E}%
U_{n}(1,h,\mathbf{t})}\right\vert \\
& \leq &\frac{|U_{n}(\varphi ,h,\mathbf{t})-\mathbb{E}U_{n}(\varphi ,h,%
\mathbf{t})|}{|U_{n}(1,h,\mathbf{t})|} \\
& &\hspace{1cm}+\;\frac{|\mathbb{E}U_{n}(\varphi ,h,\mathbf{t})|\cdot |U_{n}(1,h,\mathbf{t%
})-\mathbb{E}U_{n}(1,h,\mathbf{t})|}{|U_{n}(1,h,\mathbf{t})|\cdot |\mathbb{E}%
U_{n}(1,h,\mathbf{t})|} \\
& =:&\text{(I)}+\text{(I\negthinspace I)}.
\end{eqnarray*}%
From Theorem \ref{theo:bounded}, (\ref{stein1}) and $f_{X}$ bounded away
from zero on $J$ we get for some $\xi _{1},\xi _{2}>0$ and $c$ large enough
in $a_{n}=c(\log n/n)^{1/m}$, 
\begin{equation*}
\liminf_{n\rightarrow \infty }\sup_{a_{n}\leq h<b_{n}}\sup_{\mathbf{t}\in
I^{m}}|U_{n}(1,h,\mathbf{t})|=\xi _{1}>0,\quad \text{a.s.},
\end{equation*}%
and for $n$ large enough, 
\begin{equation*}
\sup_{a_{n}\leq h<b_{n}}\sup_{\mathbf{t}\in I^{m}}|\mathbb{E}U_{n}(1,h,%
\mathbf{t})|=\xi _{2}>0.
\end{equation*}%
Further, for $a_{n}^{\prime \prime }$ be either $a_{n}$ or $a_{n}^{\prime }$%
, we obtain readily from the assumptions (\ref{bounded}) or (\ref{unbounded}%
) on the envelope function that 
\begin{equation*}
\sup_{a_{n}^{\prime \prime }\leq h<b_{n}}\sup_{\varphi \in \mathcal{F}}\sup_{%
\mathbf{t}\in I^{m}}|\mathbb{E}U_{n}(\varphi ,h,\mathbf{t})|=O(1).
\end{equation*}%
Hence, we can now use Theorem \ref{theo:bounded} to handle (I\negthinspace
I), while for (I), depending on whether the class $\mathcal{F}$ satisfies (%
\ref{bounded}) or (\ref{unbounded}), we apply Theorem \ref{theo:bounded} or
Theorem \ref{theo:unbounded} respectively. Taking everything together we
conclude that for $c$ large enough and some $C^{\prime \prime }>0$, with
probability one, 
\begin{equation*}
\hspace{-2cm}\limsup_{n\rightarrow \infty }\sup_{a_{n}^{\prime \prime }\leq
h<b_{n}}\sup_{\varphi \in \mathcal{F}}\sup_{\mathbf{t}\in I^{m}}\frac{\sqrt{%
nh^{m}}|\hat{m}_{n,\varphi }(\mathbf{t},h)-\widehat{\mathbb{E}}\hat{m}%
_{n,\varphi }(\mathbf{t},h)|}{\sqrt{|\log h|\vee \log \log n}}\vspace{-3mm}
\end{equation*}%
\begin{eqnarray*}
\hspace{1cm}& \leq &\limsup_{n\rightarrow \infty }\sup_{a_{n}^{\prime \prime
}\leq h<b_{n}}\sup_{\varphi \in \mathcal{F}}\sup_{\mathbf{t}\in I^{m}}\frac{%
\sqrt{nh^{m}}\text{(I)}}{\sqrt{|\log h|\vee \log \log n}} \\
& &\hspace{1cm}+\;\limsup_{n\rightarrow \infty }\sup_{a_{n}^{\prime \prime
}\leq h<b_{n}}\sup_{\varphi \in \mathcal{F}}\sup_{\mathbf{t}\in I^{m}}\frac{%
\sqrt{nh^{m}}\text{(I\negthinspace I)}}{\sqrt{|\log h|\vee \log \log n}}\\
&\leq &C^{\prime \prime },
\end{eqnarray*}%
proving the assertion of Theorem \ref{theo:conv-rates}.

\appendix

\section{Appendix}

The first result below is stated as Theorem 4 in Gin\'{e} and Mason \cite%
{GM07b}, and is essentially a consequence of a martingale inequality due to
Brown \cite{Brown}. The second Theorem is a generalization of Bochner's lemma.

\begin{theor}[Theorem 4 of Gin\'{e} and Mason, 2007b]
\label{Brown} Let $X_{1},X_{2},\ldots$ be i.i.d. $S$-valued with probability
law $P$. Let $\mathcal{H}$ be a $P$--separable collection of measurable
functions $f:S^{k}\rightarrow \mathbb{R}$ and assume that $\mathcal{H}$ is $%
P $--canonical (which means that every $f$ in $\mathcal{H}$ is $P$%
--canonical). Further assume that $\mathbb{E}\Vert
f(X_{1},\dots,X_{k})\Vert_{\mathcal{H}}^{r}<\infty$ for some $r>1$, and let $%
s$ be the conjugate of $r$. Then, with $S_{n}$ defined as 
\begin{equation*}
S_{n}=\sup_{f\in\mathcal{H}}\Big|\sum_{\mathbf{i}\in
I_{n}^{k}}f(X_{i_{1}},\dots,X_{i_{k}})\Big|,\quad n\geq k,
\end{equation*}
we have for all $x>0$ and $0<c<1$, 
\begin{equation*}  \label{max}
\mathbb{P}\Big\{\max_{k\leq m\leq n}S_{m}>x\Big\}\leq \frac{\mathbb{P}%
\{S_{n}>cx\}^{1/s}(\mathbb{E}S_{n}^{r})^{1/r}}{x(1-c)}.
\end{equation*}
\end{theor}

\begin{theor}
\label{stein} Let $I=[a,b]$ be a compact interval. Suppose that $\mathcal{H}$
is a uniformly equicontinuous family of real valued functions $\varphi$ on $%
J=[a-\eta,b+\eta]^{d}$ for some $d\geq1$ and $\eta>0$. Further assume that $%
K $ is an $L_{1}$--kernel with support in $[-1/2,1/2]^{d}$ satisfying $\int_{%
\mathbb{R}^{d}}K(\mathbf{u})d\mathbf{u}=1$. Then uniformly in $\varphi\in\mathcal{H}$ and for
any sequence of positive constants $b_{n}\to0$, 
\begin{equation*}
\sup_{0<h<b_{n}}\sup_{\mathbf{z}\in I^{d}}|\varphi\ast K_{h}(\mathbf{z})-\varphi
(\mathbf{z})|\longrightarrow0,\quad\text{as }n\to\infty,
\end{equation*}
where $K_{h}(\mathbf{z})=h^{-d}K\left( \mathbf{z}/h\right) $ and 
\begin{equation*}
\varphi\ast K_{h}(\mathbf{z}):=h^{-d}\int_{\mathbb{R}^{d}}\varphi(\mathbf{x})K\left( \frac {\mathbf{z}-\mathbf{x}}{h}\right) d\mathbf{x}.
\end{equation*}
\end{theor}

\subsection{ Moment bounds}

\begin{theor}[Proposition 1 of Einmahl and Mason, 2005]
\label{em} Let $\mathcal{G}$ be a pointwise measurable class of bounded
functions with envelope function $G$ such that for some constants $%
C,\nu\geq1 $ and $0<\sigma\leq\beta$, the following conditions hold: 
\begin{itemize}
\item[(i)]$\mathbb{E} G^{2}(X)\leq\beta^{2}$; 
\item[(ii)]$\mathcal{N}(\epsilon,\mathcal{G})\leq
C\epsilon^{-\nu},\quad0<\epsilon<1$; 
\item[(iii)]$\sigma_{0}^{2}:=\sup_{g\in\mathcal{G}}\mathbb{E} g^{2}(X)\leq
\sigma^{2}$; 
\item[(iv)]$\sup_{g\in\mathcal{G}}\| g\|_{\infty}\leq\frac{1}{4\sqrt{\nu}}\sqrt{%
n\sigma^{2}/\log(C_{1}\beta/\sigma)}$,  where  $C_{1}=C^{1/\nu }\vee
e$.
\end{itemize}
Then we have for some absolute constant $A$, 
\begin{equation*}
\mathbb{E} \|\sum_{i=1}^{n}\varepsilon_{i}g(X_{i})\|_{\mathcal{G}}\leq A%
\sqrt{\nu n\sigma^{2}\log(C_{1}\beta/\sigma)},  \label{EE}
\end{equation*}
where $\varepsilon_{1},\ldots, \varepsilon_{n}$ are i.i.d Rademacher
variables independent of $X_{1},\ldots, X_{n}$.
\end{theor}

\begin{theor}[Corollary 1 of Gin\'e and Mason, 2007b]
\label{gm2007b} Let $\mathcal{F}$ be a collection of measurable functions $%
f: S^{m}\to\mathbb{R}$, symmetric in their entries with absolute values
bounded by $M>0$, and let $P$ be any probability measure on $(S,\mathcal{S})$
(with $X_{i}$ i.i.d--$P$). Assume that $\mathcal{F}$ is of VC--type with
envelope function $F\equiv M$ and with characteristics $A$ and $v$. Then for
every $m\in\mathbb{N}$, $A\geq e^{m}, v\geq1$ there exist constants $%
C_{1}:=C_{1}(m,A,v,M)$ and $C_{2}=C_{2}(m,A,v,M)$ such that for $k=1,\dots,m$%
, 
\begin{equation*}
n^{k}\mathbb{E} \|U_{n}^{(k)}(\pi_{k}f)\|_{\mathcal{F}}^{2}\leq
C_{1}^{2}2^{k}\sigma^{2}\left( \log\frac{A}{\sigma}\right) ^{k},
\end{equation*}
assuming $n\sigma^{2}\geq C_{2}\log( A/\sigma)$, where $\sigma^{2}$ is any number satisfying 
\begin{equation*}
\|P^{m}f^{2}\|_{\mathcal{F}}\leq\sigma^{2}\leq M^{2}.
\end{equation*}
\end{theor}

\subsection{Exponential inequalities}

\begin{theor}[Talagrand, 1994]
\label{tala} Let $\mathcal{G}$ be a pointwise measurable class of functions
satisfying 
\begin{equation*}
\Vert g\Vert_{\infty}\leq M<\infty,\quad g\in\mathcal{G}.
\end{equation*}
Then we have for all $t>0$,
\begin{equation*}
\hspace{-2cm}\mathbb{P}\left\{ \max_{1\leq m\leq n}\Vert\sqrt{m}\alpha_{m}\Vert _{%
\mathcal{G}}\geq A_{1}(\mathbb{E}\big\Vert \sum_{i=1}^{n}\varepsilon
_{i}g(X_{i})\big\Vert _{\mathcal{G}}+t)\right\}
\end{equation*}%
\begin{equation*}
\hspace{2cm}\leq \;2\left\{ \exp\left( -\frac{A_{2}t^{2}}{n\sigma_{\mathcal{G}}^{2}}\right)
+\exp\left( -\frac{A_{2}t}{M}\right) \right\} ,
\end{equation*}
where $\sigma_{\mathcal{G}}^{2}=\sup_{g\in\mathcal{G}}\text{Var}\left(
g(X)\right) $ and $A_{1},A_{2}$ are universal constants.
\end{theor}

We now state the exponential inequality that will permit us to control the
probability term in (\ref{ineq}), and which is stated as Theorem 5.3.14 in
de la Pe\~{n}a and Gin\'e \cite{dlPG99}.

\begin{theor}[Theorem 5.3.14 of de la Pen\~{a} and Gin\'{e}, 1999]
\label{pena} Let $\mathcal{H}$ be a $VC$--subgraph class of uniformly
bounded measurable real valued kernels $H$ on $(S^{m},\mathcal{S}^{m})$,
symmetric in their entries. Then for each $1\leq k\leq m$ there exist
constants $c_{k},d_{k}\in\;]0,\infty[$ such that, for all $n\geq m$ and $t>0$%
, 
\begin{equation*}
\Big\{\Vert n^{k/2}U_{n}^{(k)}(\pi _{k}H)\Vert _{\mathcal{H}}>t\Big\}\leq
c_{k}\exp \{-d_{k}t^{2/k}\}.
\end{equation*}%
\medskip 
\end{theor}

\noindent \textbf{Acknowledgements. } We would like to thank Uwe Einmahl for
pointing out an oversight and for very helpful suggestions. David Mason's research was partially supported by an NSF Grant. Julia Dony's research is financed by a PhD grant from the Institute for the Promotion of Innovation through Science and Technology in Flanders (IWT Vlaanderen).

\end{document}